\documentclass[11pt,oneside,reqno]{amsart}

\usepackage{amsbsy,amssymb,amscd,amsfonts,latexsym,amstext,delarray,amsmath,color,caption}
\usepackage{hyperref}
\usepackage{tikz}
\usetikzlibrary{matrix,arrows}
\usepackage{graphicx}
\usepackage{hyphenat}

\usepackage{soul}
\usepackage{amsaddr}

\usepackage{graphicx,relsize}

\usepackage{cite}

\usepackage{mathrsfs}
\usepackage[errorshow]{tracefnt}
\usepackage{fancybox}
\usepackage{amsfonts}

\usepackage{mathdots}
\usepackage{multirow}
\usepackage{array}
\usepackage{mathtools}
\usepackage{soul}
\usepackage{pict2e}
\usepackage{mathabx}
\DeclareMathAlphabet{\mathpzc}{OT1}{pzc}{m}{it}

\usepackage{float}

\let\latexcirc=\circ
\newcommand{\ccirc}{\mathbin{\mathchoice
  {\xcirc\scriptstyle}
  {\xcirc\scriptstyle}
  {\xcirc\scriptscriptstyle}
  {\xcirc\scriptscriptstyle}
}}
\newcommand{\xcirc}[1]{\vcenter{\hbox{$#1\latexcirc$}}}
\let\circ\ccirc

\makeatletter
\newcommand{\thickhline}{%
    \noalign {\ifnum 0=`}\fi \hrule height 1pt
    \futurelet \reserved@a \@xhline
}
\newcolumntype{'}{@{\hskip\tabcolsep\vrule width 1pt\hskip\tabcolsep}}
\makeatother
\newcolumntype{"}{@{\hskip\tabcolsep\vrule width 1.5pt\hskip\tabcolsep}}
\makeatother

\newcommand{\scr}{\mathscr}

\def\fg{\mathfrak{g}}
\def\fgp{\mathfrak{g}'}

\def\ie{{\it i.e.}}
\def\eg{{\it e.g.}}

\def\Z{\mathbb{Z}}
\def\boxit#1{\vbox{\hrule\hbox{\vrule\kern3pt
             \vbox{\kern3pt#1\kern3pt}\kern3pt\vrule}\hrule}}

\newcommand{\beq}{\begin{equation}}
\newcommand{\beqn}{\begin{equation*}}
\newcommand{\eeq}{\end{equation}}
\newcommand{\eeqn}{\end{equation*}}
\newcommand{\beqa}{\begin{eqnarray}}
\newcommand{\beqan}{\begin{eqnarray*}}
\newcommand{\eeqa}{\end{eqnarray}}
\newcommand{\eeqan}{\end{eqnarray*}}
\newcommand{\bdm}{\begin{displaymath}}
\newcommand{\edm}{\end{displaymath}}

\newcommand{\ba}{\begin{array}}
\newcommand{\ea}{\end{array}}

\newcommand\nn{\nonumber}

\newcommand\benu{\begin{enumerate}}
\newcommand\eenu{\end{enumerate}}
\newcommand\bit{\begin{itemize}}
\newcommand\eit{\end{itemize}}

\newtheorem{theorem}{Theorem}[section]
\newtheorem{lemma}[theorem]{Lemma}
\newtheorem{proposition}[theorem]{Proposition}

\def\Pf{\vspace{-8pt}\noindent \textbf{Proof. }}
\def\End{\mathrm{End\,}}

\def\der'{\mathfrak{der}'\,}
\def\der{\mathfrak{der}\,}
\def\str'{\mathfrak{str}'\,}
\def\str{\mathfrak{str}\,}

\def\g{\mathfrak{g}}

\def\sl{\mathfrak{sl}}
\def\gl{\mathfrak{gl}}

\def\qed{\hspace{\stretch{1}} $\Box$}

\newcommand{\de}{\delta}

\newcommand{\ad}{\text{ad}\,}

\newcommand{\adiagram}{
\begin{picture}(250,30)(45,-10)
\put(48,-10){${\scriptstyle{0}}$}
\put(88,-10){${\scriptstyle{1}}$}
\put(128,-10){${\scriptstyle{2}}$}
\put(242,-10){${\scriptstyle{n-2}}$}
\put(282,-10){${\scriptstyle{n-1}}$}
\thicklines
\put(50,10){\line(1,1){3.5}}
\put(50,10){\line(-1,1){3.5}}
\put(50,10){\line(1,-1){3.5}}
\put(50,10){\line(-1,-1){3.5}}
\multiput(50,10)(40,0){3}{\circle{10}}
\multiput(250,10)(40,0){2}{\circle{10}}
\multiput(55,10)(40,0){2}{\line(1,0){30}}
\put(135,10){\line(1,0){20}}
\put(165,10){\line(1,0){10}}
\put(185,10){\line(1,0){10}}
\put(205,10){\line(1,0){10}}
\put(225,10){\line(1,0){20}}
\put(255,10){\line(1,0){30}}
\end{picture}
}

\newcommand{\ddiagram}{
\begin{picture}(290,70)(45,-10)
\put(48,-10){${\scriptstyle{0}}$}
\put(88,-10){${\scriptstyle{1}}$}
\put(128,-10){${\scriptstyle{2}}$}
\put(242,-10){${\scriptstyle{r-3}}$}
\put(282,-10){${\scriptstyle{r-2}}$}
\put(322,-10){${\scriptstyle{r-1}}$}
\put(303,48){${\scriptstyle{r}}$}
\thicklines
\put(50,10){\line(1,1){3.5}}
\put(50,10){\line(-1,1){3.5}}
\put(50,10){\line(1,-1){3.5}}
\put(50,10){\line(-1,-1){3.5}}
\multiput(50,10)(40,0){3}{\circle{10}}
\multiput(250,10)(40,0){3}{\circle{10}}
\multiput(55,10)(40,0){2}{\line(1,0){30}}
\put(135,10){\line(1,0){20}}
\put(165,10){\line(1,0){10}}
\put(185,10){\line(1,0){10}}
\put(205,10){\line(1,0){10}}
\put(225,10){\line(1,0){20}}
\put(255,10){\line(1,0){30}}
\put(295,10){\line(1,0){30}}
\put(290,50){\circle{10}}
\put(290,15){\line(0,1){30}}
\end{picture}
}

\newcommand{\ediagram}{
\begin{picture}(330,70)(45,-10)
\put(48,-10){${\scriptstyle{0}}$}
\put(88,-10){${\scriptstyle{1}}$}
\put(128,-10){${\scriptstyle{2}}$}
\put(242,-10){${\scriptstyle{r-4}}$}
\put(282,-10){${\scriptstyle{r-3}}$}
\put(322,-10){${\scriptstyle{r-2}}$}
\put(362,-10){${\scriptstyle{r-1}}$}
\put(303,48){${\scriptstyle{r}}$}
\thicklines
\put(50,10){\line(1,1){3.5}}
\put(50,10){\line(-1,1){3.5}}
\put(50,10){\line(1,-1){3.5}}
\put(50,10){\line(-1,-1){3.5}}
\multiput(50,10)(40,0){3}{\circle{10}}
\multiput(250,10)(40,0){4}{\circle{10}}
\multiput(55,10)(40,0){2}{\line(1,0){30}}
\put(135,10){\line(1,0){20}}
\put(165,10){\line(1,0){10}}
\put(185,10){\line(1,0){10}}
\put(205,10){\line(1,0){10}}
\put(225,10){\line(1,0){20}}
\put(255,10){\line(1,0){30}}
\put(295,10){\line(1,0){30}}
\put(335,10){\line(1,0){30}}
\put(290,50){\circle{10}}
\put(290,15){\line(0,1){30}}
\end{picture}
}

\RequirePackage{hyperref}

\numberwithin{equation}{section}

\begin{document}
\noindent

\pagestyle{plain}

\vspace*{.5cm}

\frenchspacing

\title{Generators and relations for Lie superalgebras \\of Cartan type}

\author{Lisa Carbone${}^1$, Martin Cederwall${}^2$
  and Jakob Palmkvist${}^{2}$}
\address{${}^1$ Department of Mathematics, Rutgers University,\\
  110 Frelinghuysen Rd, Piscataway, NJ 05584, USA}
\address{${}^{2}$ Department of Physics, Chalmers University of Technology,\\
412 96 Gothenburg, Sweden}

\begin{abstract}   We give an analog of a Chevalley--Serre
  presentation for the Lie superalgebras $W(n)$ and $S(n)$ of Cartan
  type. These are part of a wider class of Lie superalgebras, the so-called tensor hierarchy algebras, denoted
  $W(\mathfrak{g})$ and $S(\mathfrak{g})$, where $\mathfrak{g}$
  denotes the Kac--Moody algebra 
  $A_r$, $D_r$ or $E_r$. Then $W(A_{n-1})$ and $S(A_{n-1})$ are the Lie
  superalgebras $W(n)$ and $S(n)$. The algebras $W(\mathfrak{g})$ and
  $S(\mathfrak{g})$ 
  are 
constructed from the Dynkin diagram of  
the Borcherds--Kac--Moody superalgebras $\scr B(\mathfrak{g})$ obtained by
adding a single grey node (representing an odd null root) 
to the Dynkin diagram of $\mathfrak{g}$.
  We redefine the algebras $W(A_r)$ and
$S(A_r)$ in terms of Chevalley generators and defining relations. We
  prove that all relations follow from the defining ones at level
  $\geq-2$. The analogous definitions of the algebras in the $D$- and
  $E$-series are given. In the latter case the full set of defining
  relations is conjectured.

\end{abstract}

\maketitle
\thispagestyle{empty}
\newpage

\tableofcontents

\newpage

\section{Introduction\label{IntroSection}}   

In the classification of complex finite-dimensional simple Lie
superalgebras, the classical ones are separated from those of Cartan
type, and further divided into the two classes of basic and strange
Lie superalgebras \cite{Kac77B}.  The simplest example of a Lie
superalgebra of Cartan type (from which the other ones can be
obtained) is $W(n)$, consisting of all derivations of the Grassmann
algebra with $n$ generators.  The basic Lie superalgebras are similar
to finite-dimensional simple Lie algebras, in the sense that they can
be constructed from Dynkin diagrams, where each node represents a
simple root with corresponding Chevalley generators ``$e$'' and ``$f$''. One example is $A(n-1,0)=\mathfrak{sl}(1|n)$.
In the present paper, we will show how these two representatives of two important
classes of finite-dimensional Lie superalgebras (Cartan type and classical) are in fact related to each other.
As one of our main results (Theorem~\ref{maintheorem}) we will show that
$W(n)$ can be constructed from the same Dynkin diagram as $A(n-1,0)$, but with extended sets of generators (\ref{setofgenerators})
and relations (\ref{eigen})--(\ref{eifjf0a}), (\ref{IdealJ2}). The relationship
extends to infinite-dimensional Lie superalgebras, with so-called tensor hierarchy algebras on the one hand side, and Borcherds--Kac--Moody
(BKM) superalgebras
on the other. 

Any (possibly finite\hyp{}dimensional) Kac--Moody algebra $\g$ can be
extended to a contragredient Lie superalgebra $\scr B$ by adding a
node to the Dynkin diagram such that the corresponding Chevalley
generators are  odd elements \cite{Kac77B}. When
the  simple root represented by the new node is a null root, the
node is said to be grey, and drawn as $\otimes$.  We are particularly
interested in the case where $\g$ belongs to the $A$-, $D$- or
$E$-series of Kac--Moody algebras, say $\g=X_r$ with $X$ being either
$A$, $D$ or $E$, and the grey node replaces the usual white node added
to the Dynkin diagram of $X_r$ in the extension to $X_{r+1}$. This
allows us to identify the grey node uniquely, and write $\scr B =\scr
B(\mathfrak{g})$. The Cartan matrices of $\scr B$ 
and $X_{r+1}$ only differ by the diagonal entry corresponding to the additional node.
In this case $\scr B$ is a BKM
superalgebra \cite{Borcherds,Ray}.

In \cite{Palmkvist:2013vya} a Lie superalgebra 
closely related to $\scr B$, but not of BKM type,
was defined from the same
Dynkin diagram as $\scr B$ (for finite-dimensional $\g$). 
It was called {\it tensor hierarchy algebra},
and here we denote it by $S(\mathfrak{g})$.
We also introduce a
third Lie superalgebra $W(\mathfrak{g})$ associated to the same Dynkin
diagram.
Our aim is to give a unified construction of these algebras.
To this end, Chevalley-type generators for the tensor hierarchy algebras
are defined from the same Dynkin diagram as
$\scr B(\g)$, but with an asymmetry between generators ``$e$'' and
``$f$'' in the absence of a Cartan involution.
Throughout the article, we will work
over a field $\mathbb{K}$, which can be either the complex or real
numbers.

We will focus mainly on the class of algebras obtained by taking
$X_r=A_r$, and show that $W(A_{n-1})$ and $S(A_{n-1})$ are the well known
finite\hyp{}dimensional Lie superalgebras of Cartan type, 
denoted by $W(n)$ and $S(n)$, respectively.
We will
obtain presentations for $W(A_r)$ and $S(A_r)$ by giving an analog of Chevalley
generators and defining relations.
We will also comment on the cases $X_r=D_r$ and $X_r=E_r$.
The latter case is interesting, and provides the main motivation from
mathematical 
physics for this investigation, due to a deep, and to a large extent
unexplored, connection to
generalised and extended geometry
(see
\cite{Palmkvist:2011vz,deWit:2008ta,Cederwall:2015oua,Berman:2012vc,Cederwall:2013naa,
Palmkvist:2012nc,Bossard:2017wxl,BossardEtAl:2017,CederwallPalmkvist:2017}).
The level decompositions of $W(E_r)$  and $S(E_r)$ correctly
predict, for example,  the embedding tensor used in the construction of gauged
supergravities. Given the geometric character of the classical
definition of $W(n)$ as operating on forms, one may envisage a similar
r\^ole for $W(E_r)$ in exceptional geometry. Hopefully, the tensor
hierarchy algebra can be shown to provide an underlying structure, on
which the concept of exceptional geometry depends. The relationship between tensor hierarchies and Leibniz
algebras recently studied in \cite{Lavau:2017tvi} might be useful in this respect.
From a mathematical
perspective, we expect that the present construction will be useful
\eg\ for addressing
questions concerning the combinatorics and representation theory of Cartan type
superalgebras. 

Sections \ref{IntroSection} and \ref{sec:graded} contain a review of
the BKM and Cartan type
superalgebras used in the paper, and introduce some notation.
In Section \ref{WgSection}, we define the Lie superalgebras
$\widetilde W(\fg)$ and $\widetilde S(\fg)$ in terms of Chevalley--Serre-like
generators and relations, from which $W(\fg)$ and $S(\fg)$
are obtained by factoring out the maximal ideal intersecting the local part trivially.
In Section \ref{IdealJSection}, we construct this maximal
ideal for the case $\fg=A_{n-1}$. We end in Section 5 with a discussion of and a
conjecture for the $D$
and $E$ cases in Section \ref{DESection}, where  the
identification $S(D_r)\simeq H(2r)$ is made. Details about the root system of $W(A_{n-1})=W(n)$ are given in an appendix.

\subsection*{Acknowledgments} JP would like to thank Axel Kleinschmidt for valuable discussions,
and the CERN theory division for providing a stimulating environment.
MC and JP would like to thank Rutgers University for hospitality, and
LC and JP would like to thank Institut des Hautes \'Etudes Scientifiques (IH\'ES) for its generous support.
The work of MC and JP is
supported by the Swedish Research Council, project no. 2015-04268.
The work of LC is supported by  Simons Foundation Collaboration Grant, no. 422182.
Part of the work of JP was done at the Mitchell Institute for Fundamental Physics and Astronomy at 
Texas A\&M University, and supported in part by NSF grant PHY-1214344.

\section{\texorpdfstring{$\mathbb{Z}$}{Z}-graded Lie superalgebras}
\label{sec:graded}

In this section we review some basic definitions and results from
Section 1.2 in \cite{Kac77B}. We refer to this article for further details. 

First we recall that a Lie superalgebra is a $\mathbb{Z}_2$-graded
vector space $G=G_{(0)} \oplus G_{(1)}$ with a bilinear bracket 
\begin{align}
G \times G &\to G\,, & (x,y) &\mapsto [x,y] 
\end{align}
satisfying $[G_{(i)},G_{(j)}] \subseteq G_{((i+j)\text{ mod } 2)}$, and the identities 
\begin{align}
[x,y]&=-(-1)^{|x||y|}[y,x]\,, \label{supersymmetry}\\
[x,[y,z]]&=[[x,y],z]-(-1)^{|x||y|}[y[x,z]]\,, \label{jacobi-identity}
\end{align}
where $|x|=0$ if $x\in G_{(0)}$ and $|x|=1$ if $x\in G_{(1)}$.

A $\mathbb{Z}$-{\it grading} of the Lie superalgebra $G$ is a
decomposition of $G$  into a direct sum of subspaces $G_i$ for all
integers $i$,  called {\it levels}, such that $[G_i,G_j]\subseteq
G_{i+j}$.  In all cases we coinsider in this paper, we have $[G_i,G_j]= G_{i+j}$.

Whenever we use the notation $G_i$ for subspaces of an algebra $G$, we assume
 a $\mathbb{Z}$-grading of $G$. We will also use the notation
$G_\pm = \bigoplus_{i \in \mathbb{Z}_\pm}G_i$.  The
$\mathbb{Z}$-grading is said to be {\it consistent} if $G_i \subseteq G_{(i\text{ mod } 2)}$.

It follows from the relations $[G_i,G_j]\subseteq G_{i+j}$ that the subspace
$G_0$ of any $\mathbb{Z}$-graded Lie superalgebra $G$ is a subalgebra,
which is a Lie algebra if the $\mathbb{Z}$-grading is consistent, and 
all subspaces $G_i$ can be considered as $G_0$-modules.

A $\mathbb{Z}$-graded Lie superalgebra can be constructed
from a $\mathbb{Z}_2$-graded vector space \linebreak ${g}=g_{(0)} \oplus g_{(1)}$
with a consistent decomposition into a direct sum of subspaces\linebreak 
${g}=g_{-1} \oplus g_{0} \oplus g_1$ and a bilinear bracket defined
for all pairs of elements in $G$ such that not both of them have
nonzero components in the same subspace $g_1$ or $g_{-1}$. If
$[g_i,g_j]\subseteq g_{i+j}$ 
and the identities (\ref{supersymmetry})--(\ref{jacobi-identity}) 
are satisfied whenever the brackets are defined, then ${g}$ is
a {\it local Lie superalgebra}. 

Clearly, any $\mathbb{Z}$-graded Lie superalgebra $G=\bigoplus_{i \in
  \mathbb{Z}}G_i$ gives rise to a local Lie superalgebra $G_{-1}
\oplus G_0 \oplus G_1$, which is called the {\it local part} of $G$.
If a subspace $g'$ of a local Lie
superalgebra $g$ itself is a local Lie superalgebra with respect to
the bracket and the decomposition $g'=g'{}_{-1} \oplus g'{}_{0} \oplus
g'{}_{1}$ inherited from $g$, then we call $g'$ a {\it local
  subalgebra} of the local Lie superalgebra $g$.

Given a local Lie superalgebra $g=g_{-1}\oplus g_0\oplus g_{1}$ we can
construct {\it maximal} and {\it minimal} Lie superalgebras with $g$
as the local part. The maximal one is defined as $\widetilde G = F(g)/I$
where $F(g)$ is the free Lie superalgebra generated by $g$, and $I$
the ideal generated by the commutation relations in $g$.  The maximal Lie superalgebra can then be
decomposed as $\widetilde G = \widetilde G_-  \oplus\   G_0\   \oplus  \widetilde G_+$,
where $\widetilde G_\pm$ is the free Lie superalgebra generated by $\widetilde
G_{\pm1}=g_{\pm1}$. The minimal Lie superalgebra with  local part $g$
is defined as $G=\widetilde G / J$, where $J$ is the maximal homogeneous
ideal of $\widetilde G$ intersecting the local part trivially. Thus
$G_{\pm1}=\widetilde G_{\pm1}=g_{\pm1}$.  Maximality (respectively minimality) here means that
any isomorphism between the local parts of $\widetilde G$ (respectively $G$) and any
other $\mathbb{Z}$-graded Lie superalgebra $G'$ can be extended to a
surjective homomorphism $\widetilde G\to G'$ (respectively $G\to G'$)
\cite{Kac68,Kac77B}.

\subsection{The local Lie superalgebra associated to a vector space} \label{universal}

Any $\mathbb{Z}_2$-graded vector space $U_1$ gives rise to a local Lie superalgebra in the following way.
Set
\begin{align}
U_0&=\End U_1\,,& U_{-1}&={\rm Hom}(U_1,U_0)\,.
\end{align}
Then $U_1 \oplus U_0 \oplus U_{-1}$ 
is a local Lie superalgebra, denoted $u (U_1)$, with the bracket given by the following relations,
\begin{align} \label{local1}
[x_0,y_1]&=x_0 (y_1)\,, & 
[x_0,y_0]&=x_0\ccirc y_0-(-1)^{|x||y|}y_0\ccirc x_0\,,\nn\\
[x_{-1},y_1]&=x_{-1}(y_1)\,, & 
[x_0,y_{-1}]&=x_0\ccirc  y_{-1} - (-1)^{|x||y|} y_{-1}\ccirc  x_0\,,
\end{align}
where $x_i$ and $y_i$ belong to $U_i$ ($i=0,\pm1$), and have $\mathbb{Z}_2$-degrees $|x|$ and $|y|$, respectively. In particular, this means
that $U_0$, as a Lie superalgebra, is $\gl(U_1)$.
Let $K_a$ be a basis of $U_1$, for $a=1,2,\ldots,{\rm dim }\,U_1$. 
Then we have bases $K^a{}_b$ of $U_0$ (with $a,b$ indexing rows and
columns of a matrix in $\End U_1$) and $K^{a,b}{}_c$  of $U_{-1}$ (with $a$ indexing the basis of the dual $U_1^*$ and $b,c$ indexing the basis of $U_0$). 
When $U_1$ is an odd vector space, so that $U_0=\mathfrak{gl}(U_1)$ is a Lie algebra, the commutation relations for the basis elements
that follow from (\ref{local1}) are
\begin{align} \label{unilocalcommrel}
[K_a,K^b{}_c]&= \de_a{}^b K_c\,, & [K^a{}_b,K^c{}_d]&=\de_b{}^c
K^a{}_d - \de_d{}^a K^c{}_b\,, \nn\\ [K_a,K^{b,c}{}_d]&=\de_a{}^b
K^c{}_d\,, & [K^a{}_b,K^{c,d}{}_e]&=\de_b{}^c K^{a,d}{}_e+\de_b{}^d
K^{c,a}{}_e - \de_e{}^a K^{c,d}{}_b\,,
\end{align}
where $\de_a{}^b$ is a Kronecker delta.

From this local Lie superalgebra, and its local subalgebras, we can
construct the maximal and minimal Lie superalgebras as above. As we
    will see, the BKM algebras described in the next section,
    and the associated
    tensor hierarchy algebras, 
    can be constructed in this way.

\subsection{Borcherds--Kac--Moody superalgebras}

We are particularly interested in
Lie superalgebras $\scr B$ that are minimal
with respect to a consistent $\mathbb{Z}$-grading and with local part $\scr B_{-1} \oplus \scr B_0 \oplus
\scr B_1$, where the subalgebra $\scr B_0$ is the direct sum of a simple and
simply laced (possibly finite-dimensional) Kac--Moody algebra
$\mathfrak{g}$ and a one-dimensional center of $\scr B_0$, and the
representation of $\mathfrak{g}$ on $\scr B_1$ is fundamental and dual
to the one on $\scr B_{-1}$. Then in the Dynkin diagram of $\scr B$, 
there is a grey
node connected to one of the white nodes in the Dynkin
diagram of $\mathfrak{g}$ by a single line.

Let $r$ be the rank of $\mathfrak{g}$. We use a labelling $1,2,\ldots,r$ of the nodes in the Dynkin diagram
of $\mathfrak{g}$ such that
the grey node is connected to node 1, and can be included as node 0 in
an extended labelling $0,1,2,\ldots,r$. Then $\scr B$ has a
corresponding Cartan matrix $B_{ab}$ ($a,b=0,1,\ldots,r$) such that
\begin{align}
B_{00}=0, \qquad B_{01}=B_{10}=-1,\nn
\end{align}
\begin{align}
B_{0i}=B_{i0}=0 \qquad 
 (i=2,\ldots,r). \label{cartan-B}
\end{align}
Assuming $B_{ab}$ to be non-degenerate, $\scr B$
can be constructed as the Lie superalgebra
generated by $2(r+1)$ elements $e_a,f_a$ (odd if $a=0$, even otherwise)
modulo the Chevalley--Serre relations
\begin{align} \label{chev-rel}
[h_a,e_b]&=B_{ab}e_b\,, & [h_a,f_b]&=-B_{ab}f_b\,, & [e_a,f_b]&=\delta_{ab}h_b\,,
\end{align}
\begin{align} \label{serre-rel}
(\text{ad }e_a)^{1-B_{ab}}(e_b)&=(\text{ad }f_a)^{1-B_{ab}}(f_b)=0
\qquad\quad(a\neq b)\,
\end{align}
($a,b=0,1,2,\ldots,r$), where the elements $h_a = [e_a,f_a]$ span an
abelian  Cartan subalgebra. Thus the Chevalley--Serre relations
for $\scr B$ take the same form as those for $\mathfrak{g}$ (keeping
in mind that the index set is extended, and that the bracket
$[e_0,f_0]$ is symmetric, but it is often relevant to consider $\scr B$
as a special case of a BKM superalgebra
with additional relations in the general case. 
We refer to \cite{Ray} for details
about general BKM superalgebras.

A comment on notation:
We write $\scr B$ as ${\scr B}(\fg)$ 
when we want to emphasise the underlying Lie algebra $\fg$.
Strictly speaking, it is
not sufficient to know the Lie algebra $\fg$ 
in order to construct $\scr B(\fg)$ from it;
the data specifying ${\scr B}(\fg)$ is a choice of $\fg$ together with
a choice of node $1$ (the node connecting to the grey node
$0$). In the series of greatest interest to us, the $A$-, $D$- and $E$-series,
our default choice of node $1$  is the node connected to the additional white node 
added in order to obtain the next algebra in the series. However, other
choices are possible. The same comment applies to the Lie superalgebras
$W(\fg)$ and $S(\fg)$ that we will later apply to the same Dynkin diagram as $\scr B(\fg)$.

Throughout the rest of the paper, the indices $a,b,\ldots$ will take the $r+1$ values $0,1,\ldots,r$,
where $r$ is the rank of the Lie algebra $\fg$. When $\fg$ belongs to the $A$-series it is
convenient to set $n=r+1$, and thus the indices $a,b,\ldots$ will take the $n$ values $0,1,\ldots,n-1$. The range of  indices $i,j,\ldots$ may vary, and will be specified explicitly whenever they appear (if not obvious from the context).

It follows from (\ref{cartan-B}) that ${\rm det}\, B= - {\rm det}\, A'$ where $A'$ is the Cartan
matrix obtained by removing the first two rows and columns
(corresponding to removing nodes 0 and 1 from the Dynkin diagram), as is evident in (\ref{cartan-A-case}) below.
The subalgebra of $\fg$ with the Cartan matrix $A'$ obtained in this way will play an important role later,
and we denote it by $\fg'$. Since we assume $B$ to be non-degenerate, this means that $\fg'$ is simple.

As is the case for  Kac--Moody algebras, the nodes in the Dynkin diagram of $\scr
B$ correspond not only to generators $e_a$ and $f_a$, but also to simple
roots $\alpha_a$ that form a basis of the dual space of 
the Cartan subalgebra, and the Cartan matrix defines an inner product
on this space up to an overall normalisation.  Since we only
consider simply laced Dynkin diagrams, we do not have to symmetrise
the Cartan matrix, and we fix the normalisation as follows:
\begin{align}
(\alpha_a,\alpha_b)=B_{ab}.
\end{align}
Thus $(\alpha_0,\alpha_0)=0$, so $\alpha_0$ is a null root, whereas
$(\alpha_i,\alpha_i)=2$ for $i=1,2,\ldots,r$.

The inner
product also defines the  basis $\Lambda_1$, $\Lambda_2$,$\ldots,\Lambda_r$ of
fundamental weights of $\fg$, dual to the basis
$\alpha_1,\alpha_2,\ldots,\alpha_r$ of simple roots, with
$(\Lambda_i,\alpha_j)=\delta_{ij}$.

In the consistent $\mathbb{Z}$-grading of $\scr B$ mentioned above,
the odd generators $e_0$ and $f_0$ belong to $\scr B_{1}$ and $\scr
B_{-1}$, respectively, and  the even generators belong to $\scr B_{0}$.
As a $\mathfrak{g}$-module, $\scr B_1$ has a lowest weight vector
$e_0$ with weight $-\Lambda_1$, and $\scr B_{-1}$ has a highest weight
vector $f_0$ with weight $\Lambda_1$. We will often use the following notation for
weights where $\Lambda_1=(10\ldots 0)$, and  the entries are
the Dynkin labels, \ie, the coefficients of the weight in the
basis $\Lambda_1,\Lambda_2,\ldots,\Lambda_r$. 
We will occasionally
use this notation not only for
the weight itself, but also for the irreducible representation with that
highest weight. Occasionally we will instead use the notation $R(\lambda)$ for the
irreducible module with highest weight $\lambda$. Thus
\begin{align}
\scr B_{-1} = R(\Lambda_{1}) = (10\ldots 0)\,.
\end{align}

\subsection{The associated tensor hierarchy algebras} \label{THA-subsection} 
In this subsection, we recall the tensor hierarchy algebras of \cite{Palmkvist:2013vya} and we consider variations on the construction of $\scr B$ which allow us to modify the zero and negative levels of $\scr B$ without affecting the positive levels. 
We assume a $\mathbb Z$-grading of $\scr B$ as in the preceding subsection.

Consider the local Lie superalgebra 
$u(\scr B_1)=U_{-1} \oplus U_0 \oplus U_1$
associated to the subspace $\scr B_1$ of $\scr B$, as defined in Section~\ref{universal}.
Thus
\begin{align}
U_{-1}&={\rm Hom}(\scr B_1,\End \scr B_1), & U_0 &= \End\scr B_1, & U_1 &= \scr B_1.
\end{align}
The adjoint action in $\scr B$ of $\scr B_{0}=\fg \oplus \mathbb{K}$ on $\scr B_1$ provides an embedding of
$\fg \oplus \mathbb{K}$ into $\End \scr B_1$. By this embedding we can restrict $U_{-1}$
to the subspace 
\begin{align}
V_{-1}={\rm Hom}(\scr B_1,\fg\oplus\mathbb{K}),
\end{align}
and consider the local subalgebra
\begin{align}
v(\scr B_1,V_{-1})=V_{-1} \oplus V_0 \oplus V_1
\end{align}
of $u(\scr B_1)$ generated by $V_1=\scr B_1$ and $V_{-1}$. We then have 
\begin{align}
V_0 =
    [V_1,V_{-1}]= \mathfrak{g}\oplus\mathbb{K},
\end{align}
and as a $\mathfrak{g}$-module,
    $V_{-1}$ is isomorphic to $\widebar{\scr B_1} \otimes (\mathfrak{g}\oplus\mathbb{K})$, the tensor product
    of the dual of $\scr B_1$ and the adjoint of $\mathfrak{g}\oplus\mathbb{K}$.  
    
    We can restrict 
    $U_{-1}$ further by varying our choice of $V_{-1}$. For example, we may choose $V_{-1}$ to be isomorphic to a submodule of the tensor
    product $\widebar{\scr B_1} \otimes (\mathfrak{g}\oplus\mathbb{K})$, where $\widebar{\scr B_1}$ is the dual of $\scr B_1$ (and thus isomorphic to $\scr B_{-1}$).
    In particular, we may take $V_{-1}$ to be isomorphic to $\widebar{\scr B_{1}}$ as a $\mathfrak{g}$-module, noting that $\widebar{\scr B_{1}}$ has multiplicity 2 as a
    $\mathfrak{g}$-module in the tensor product.
Taking $V_{-1}$ to be a particular linear combination of the two $\widebar{\scr B_{1}}$ modules in the tensor product,
we get back the BKM superalgebra $\scr B$ of Section~\ref{universal}, that is $\scr B=V=V(\scr B_1,V_{-1})$, where $V=V(\scr B_1,V_{-1})$ is the minimal Lie superalgebra with  local part $v(\scr B_1,V_{-1})$. 
Hence there  is a $\mathbb Z$-grading of $V$ such that the $i$-component indeed equals $V_{i}$ for $i=0,\pm 1$.

For finite-dimensional $\fg$ it is also possible to choose $V_{-1}$ to be 
the direct sum of $\widebar{\scr B_{1}}$ and an additional module contained in the tensor product $\widebar{\scr B_{1}} \otimes \mathfrak{g}$
such that $V_2 \subseteq \scr B_2$. The maximal additional submodule
of $\widebar{\scr B_{1}} \otimes \mathfrak{g}$ for which this holds is 
known as the {\it embedding tensor representation} in the $\fg$-covariant formulation of gauged supergravity 
with broken global symmetry $\fg$ \cite{deWit:2008ta}.

If instead of choosing $V_0=\mathfrak{g}\oplus\mathbb{K}$, we have $V_0=\mathfrak{g}$, and we  require that
$V_{-1}$ is the maximal subspace of $U_{-1}$ such that $V_2 \subseteq \scr B_2$, then 
$V_{-1}$ consists only of the embedding tensor representation (not the direct sum with $\widebar{\scr B_{1}}$) and
$V$ is precisely what was called the {\it tensor hierarchy algebra} of $\fg$ in \cite{Palmkvist:2013vya}.
 Here we extend the definition of tensor hierarchy algebra to include the case where $V_0=\mathfrak{g}\oplus\mathbb{K}$.
In the application to gauged supergravity, the difference between the two algebras depends on whether or not the so-called 
trombone gauging is taken into account. 

In Section \ref{WgSection} we will define Lie superalgebras $W=W(\fg)$ and $S=S(\fg)$
associated to $\fg$
and show that they agree (with some minor exceptions) with the tensor hierarchy algebras of $\fg$ in the case when $\fg$ is finite-dimensional ($S$ being the original one,
and $W$ the extended version), although the definition of these algebras that we will give in Section \ref{WgSection}
is very different from the construction above (and more general since it can be applied also to infinite-dimensional $\fg$).
Thus we can write $W_0=\mathfrak{g}\oplus\mathbb{K}$ and $S_0=\mathfrak{g}$.

To understand the condition $V_2 \subseteq \scr B_2$,
we recall that the minimal Lie
superalgebra $V=V(\scr B_1,V_{-1})$ can be obtained from the maximal
one $\widetilde V=\widetilde V(\scr B_1,V_{-1})$ by factoring out the maximal ideal
that intersects the local part trivially, and that $\widetilde V_+$ is
the free Lie superalgebra generated by the odd subspace $V_1=\scr
B_1$. Thus $\widetilde V_2$ is the full symmetric tensor product of two $\scr
B_1$ modules, and decomposes into a direct sum of $\scr B_2$ and
another submodule, which we denote by $\scr B_2{}^{\mathsf{c}}$ (unless $\scr B_2=0$,
in which case we set $\scr B_2{}^{\mathsf{c}}=\widetilde V_2=$). The
condition $V_2 \subseteq \scr B_2$ now means that $\scr
B_2{}^{\mathsf{c}}$ must be contained in the maximal ideal of $\widetilde V$
that intersects the local part trivially, and is thus equivalent to
$[V_{-1},\scr B_2{}^\mathsf{c}]=0$.

Any irreducible submodule of $U_{-1}$ that is not contained in
$V_{-1}$, and thus gives a nonzero bracket with $\scr
B_2{}^\mathsf{c}$, must be dual to a submodule of $\scr B_2{}^\mathsf{c}
\otimes \widebar{\scr B_{1}}$.  Thus the modules $W_{-1}$ and $S_{-1}$ can
be determined by decomposing $\widebar{\scr B_{1}} \otimes
(\mathfrak{g}\oplus \mathbb{K})$ and $\widebar{\scr B_{1}} \otimes
\mathfrak{g}$, respectively, into irreducible submodules, and subtracting the
overlap with $\widebar{\scr B_2{}^\mathsf{c}} \otimes \scr B_1{}$.
For $W_{-1}$, this is precisely the computation that determines the
torsion representation in exceptional geometry, 
defined as the part of the affine connection that transforms with the
generalised Lie derivative under a generalised 
diffeomorphism \cite{Cederwall:2013naa}. In this context, the $\scr
B_2$ representation is the one that appears in the so called section
condition. 

\subsection{The BKM superalgebra \texorpdfstring{$\scr B(A_{n-1}) =
    A(n-1,0)$}{B(A(n-1))=A(n-1,0)}} 
\begin{figure}
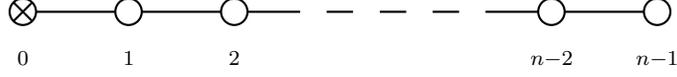

\begin{center}
\adiagram 
\end{center}
\caption{\it \label{superADynkin}The Dynkin diagram of $\scr B(A_{n-1})=A(n-1,0)$.}
\end{figure}

We now explicitly consider the case
$\mathfrak{g}=A_{n-1}=\sl(n)$, where the BKM algebra $\scr B(\fg)$
is the finite-dimensional Lie superalgebra
$\scr B=A(n-1,0)=\sl(1|n)$ with the Dynkin
diagram given in Figure \ref{superADynkin}
and with Cartan matrix  
\begin{align} \label{cartan-A-case}
B_{IJ}=\begin{pmatrix}
0 & -1 & 0 &\cdots& 0 & 0\\
-1 & 2 & -1 &\cdots& 0 & 0\\
0 & -1 & 2 &\cdots& 0 & 0\\
\vdots &\vdots &\vdots & \ddots & \vdots & \vdots\\ 
0 & 0& 0 & \cdots & 2 & -1\\
0 & 0& 0 & \cdots & -1 & 2
\end{pmatrix}\,.
\end{align}

The consistent $\mathbb{Z}$-grading of $\scr B$ is in this case a
3-grading, which means that the full Lie superalgebra coincides with its local part, 
\begin{align} \label{3gradingofA}
A(n-1,0) = \mathscr B_{-1} \oplus \mathscr B_{0} \oplus \mathscr B_{1}\,.
\end{align}
The subalgebra $\scr B_0$ is $\mathfrak{sl}(n)\oplus\mathbb{K}=\mathfrak{gl}(n)$,
and the basis elements can thus be written as $\mathfrak{gl}(n)$
tensors as in Table \ref{superAtable}.
\begin{table}[H]
\begin{align*}
\begin{array}{|c|c|c|}
\hline
\text{level} & \text{basis} & \text{$A_{n-1}$ representation}\\
\hline
1 & E_a & (00\cdots 01)\\
0 & G^a{}_b & (10\cdots 01) \oplus (00\cdots 00)\\
-1 & F^a & (10\cdots 00)\\
\hline
\end{array}
\end{align*}
\vspace{4pt}
\caption{\it \label{superAtable}The $\Z$-grading of $A(n-1,0)$.}
\vspace{-16pt}
\end{table}
The commutation relations are
\begin{align}
[G^a{}_b,G^c{}_d]&=\de_b{}^cG^a{}_d-\de_d{}^aG^c{}_b\,,  & [E_a,F^b]&=-G^b{}_a +\de_a{}^b G\,,\nn
\end{align}
\begin{align}
[G^a{}_b,F^c]&= \de_b{}^c F^a\,,& [G^a{}_b,E_c]&= -\de_c{}^a E_b\,, & [E_a,E_b]&=[F^a,F^b]=0\,,
\label{A(n-1,0)commrel}
\end{align}
where $G=\sum_{a=0}^{n-1}G^a{}_a$.
Identifying the Chevalley generators,
\begin{align} \label{chev1}
  e_0 &= E_0\,, & f_0 &= F^0\,,
  & h_0&=G^1{}_1 + G^2{}_2 + \cdots + G^{n-1}{}_{n-1} = G-G^0{}_0\,,
\end{align} 
\begin{align} \label{chev2}
e_i &= G^{i-1}{}_{i}\,, & f_i &= G^{i}{}_{i-1}\,, & h_i&=G^{i-1}{}_{i-1}-G^{i}{}_{i}\,,  & (i&=1,2,\ldots,n-1)
\end{align}
the commutation relations (\ref{A(n-1,0)commrel}) follow from the
Chevalley--Serre relations (\ref{chev-rel})--(\ref{serre-rel}). 

\subsection{The Cartan type Lie superalgebra
  \texorpdfstring{$W(n)$}{W(n)}}
    
Let $\Lambda(n)$ be the Grassmann superalgebra with generators
$\xi^0,\xi^1,\ldots,\xi^{n-1}$, that is, the associative superalgebra generated
by these elements modulo the relations $\xi^a\xi^b = -
\xi^b\xi^a$.  The Cartan type superalgebra $W(n)$ is
the derivation superalgebra of
$\Lambda(n)$, with basis elements
\begin{align}
K^{a_1 \cdots a_p}{}_b = \xi^{a_1} \cdots \xi^{a_p}
\frac{\partial}{\partial\xi^b}
\end{align}
acting on a monomial $\xi^{c_1}\cdots\xi^{c_q}$ by a contraction,
\begin{align} \label{contraction}
K^{a_1 \cdots a_p}{}_b \quad : \quad \xi^{c_1}\cdots\xi^{c_q}
\mapsto q\,\delta_b{}^{[c_1}\xi^{|a_1} \cdots \xi^{a_p|}\xi^{c_2}\cdots \xi^{c_q]}\,,
\end{align}
where the square brackets denote antisymmetrisation with total weight 1, excluding indices enclosed between vertical bars.
It is easy to see that $W(n)$ has a consistent $\mathbb{Z}$-grading
where $W_{-p+1}$ has a
basis of elements $K^{a_1 \cdots a_p}{}_b$ which are antisymmetric in
all upper indices, and thus $W_{-p+1}=0$ for $p>n$ (or $p<0$) as can
be seen in Table \ref{WnTable}. We note that negative and positive levels
are reversed compared to the usual conventions.  
\begin{table}[H]
\begin{align*}
\begin{array}{|c|c|c|c|}
\hline
\text{level} & \text{basis} & \text{$A_{n-1}$ representation} &
\text{dimension}\\ 
\hline
1 & K_a & (00\cdots 01)& n\\
0 & K^a{}_b &(10\cdots 01) \oplus (00\cdots 00)& n^2\\
-1 & K^{ab}{}_c & (010\cdots 01)\oplus (10\cdots 00)&\binom{n}{2}\cdot n\\
-2 & K^{abc}{}_d & (0010\cdots 01)\oplus (010\cdots 00)&\binom{n}{3}\cdot n\\
\vdots & \vdots & \vdots& \vdots\\
-n+2  & K^{a_1 \cdots a_{n-1}}{}_b  & (00\cdots 02)\oplus (00\cdots
010) &\binom{n}{n-1}\cdot n = n^2\\ 
-n+1  & K^{a_1 \cdots a_n}{}_b& (000\cdots 01)& n\\
\hline 
\end{array}
\end{align*}
\vspace{4pt}
\caption{\it \label{WnTable}The $\Z$-grading of $W(n)$.}
\vspace{-16pt}
\end{table}
The commutation relations are
\begin{align} \label{gencommrel}
[K^{a_1 \cdots a_p}{}_c, K^{b_1 \cdots
    b_q}{}_d]&=q\,\de_c{}^{[b_1|}K^{[a_1\cdots a_p] |b_2\cdots
    b_q]}{}_d -p\,\de_d{}^{[a_p}K^{a_1\cdots a_{p-1}] b_1 \cdots
  b_q}{}_c\,.
\end{align}

The Lie superalgebra $W(n)$ has a subalgebra $S(n)$ spanned by the traceless linear
combinations
\begin{align}
\hat K^{a_1\cdots a_p}{}_c &= K^{a_1\cdots a_p}{}_c
-\sum_{d=0}^{n-1} \frac{p}{n-p+1}\de^{[a_p}_c K^{a_1\cdots a_{p-1}]d}{}_d\,,
\end{align}
satisfying
\begin{align}
\sum_{c=0}^{n-1} \hat K^{a_1\cdots a_{p-1}c}{}_c = 0\,.
\end{align}
The commutation relations in $S(n)$ are
\begin{align}
[\hat K^{a_1\cdots a_p}{}_c,\hat K^{b_1\cdots
    b_q}{}_d]&=q\,\de_c{}^{[b_1|}\hat K^{[a_1\cdots a_p] |b_2\cdots
    b_q]}{}_d -p\,\de_d{}^{[a_p}\hat K^{a_1\cdots a_{p-1}] b_1 \cdots
  b_q}{}_c\nn\\ &\quad\,-\frac{p(q-1)}{n-p+1}\de_c^{[a_p}\hat K^{a_1
    \cdots a_{p-1}]b_1\cdots
  b_q}{}_d\nn\\ &\quad\,+\frac{q(p-1)}{n-q+1}\de_d^{[b_1|} \hat
  K^{[a_1\cdots a_p]| b_2\cdots b_{q-1}]}{}_c\,,
\end{align}
in particular
\begin{align}
[K_c,\hat K^{b_1\cdots b_q}{}_d]=q\de^{[b_1}_c \hat K^{b_2 \cdots b_q]}{}_d -
\frac{q}{n-q+1}\de^{[b_1}_d \hat K^{b_2 \cdots b_q]}{}_c\,,
\end{align}
and the $\mathbb{Z}$-grading inherited from $W(n)$ is given by Table
\ref{SnTable}.
\begin{table}[H]
\begin{align*}
\begin{array}{|c|c|c|c|}
\hline
\text{level} & \text{basis} & \text{$A_{n-1}$ representation} &
\text{dimension}\\ 
\hline
1 & K_a & (00\cdots 01)& n\\
0 & \hat K^a{}_b &(10\cdots 01) & n^2-1\\
-1 & \hat K^{ab}{}_c & (010\cdots 01)&\binom{n}{2}\cdot n-n\\
-2 & \hat K^{abc}{}_d & (0010\cdots 01)&\binom{n}{3}\cdot n-\binom{n}{2}\\
\vdots & \vdots & \vdots& \vdots\\
-n+2  & \hat K^{a_1 \cdots a_{n-1}}{}_b  & (00\cdots 02)&\frac12{n(n+1)}\\
\hline 
\end{array}
\end{align*}
\vspace{4pt}
\caption{\it \label{SnTable}The $\Z$-grading of $S(n)$.}
\vspace{-16pt}
\end{table}

\subsubsection{Towards generators and relations for
  \texorpdfstring{$W(n)$}{W(n)}} 

We note that the subalgebra of $W(n)$ generated by $K_a$, $K^{a}{}_b$
and $K^a  \equiv \sum_{b=0}^{n-1} K^{ab}{}_b$ is isomorphic to $A(n-1,0)$,
with an injective homomorphism $\psi: A(n-1,0) \to W(n)$  given by
\begin{align}
\psi(E_a)&=K_a\,, & \psi(G^a{}_b)&=K^a{}_b\,, & \psi(F^a)&=K^a\,. 
\end{align}
This follows from comparing the commutation relations
(\ref{A(n-1,0)commrel}) with those for $K_a$, $K^{a}{}_b$ and $K^a$
coming from (\ref{gencommrel}).
It is therefore natural to ask whether the set
  of Chevalley generators of $A(n-1,0)$ can be extended by generators
  corresponding to the traceless part of $K^{ab}{}_c$, and whether
  $W(n)$ then can be constructed from the extended set of generators
  by some relations extending the Chevalley--Serre relations. In order
  to investigate this we redefine the Chevalley generators
  (\ref{chev1})--(\ref{chev2}) as elements in $W(n)$ by the
  injective homomorphism $\psi$ above,
\begin{align}
  e_0 &= K_0\,, & f_0 &= K^0\,,
  & h_0&=K^1{}_1 + K^2{}_2 + \cdots + K^{n-1}{}_{n-1} = K-K^0{}_0\,,\nn
\end{align} 
\begin{align} \label{W-chev2}
e_i &= K^{i-1}{}_{i}\,, & f_i &= K^{i}{}_{i-1}\,, & h_i&=K^{i-1}{}_{i-1} -K^{i}{}_{i}\,,
\end{align}
for $i=1,\ldots,n-1$, where $K=\sum_{a=0}^{n-1}K^a{}_a$.  Accordingly, we define the Cartan subalgebra of
$W(n)$ as the subalgebra spanned by the elements
$h_0,h_1,\ldots,h_{n-1}$.  Usually the Cartan subalgebra of $W(n)$ is
defined as the Cartan subalgebra  of the $A_{n-1}$ subalgebra, which is spanned by
$h_1,\ldots,h_{n-1}$,
but from our point
of view, it is natural to include also $h_0$.  We can thus define roots
in the usual way with respect to this Cartan subalgebra.

We describe the root system of $W(n)$ in the appendix, in general, and explicitly for $n=3$ and $n=4$. As can be seen for $W(3)$ in Table \ref{w3table},
there are not only positive and negative roots, but also a root $-\alpha_0+\alpha_2$ mixing positive and negative coefficients
for the simple roots. This is a general feature of $W(n)$ root systems.
Another remarkable feature in which the Cartan type superalgebras differ from the BKM superalgebras
is that the simple root $-\alpha_0$ has multiplicity greater than one (as usual, the {multiplicity} of a root is the dimension of its root space).
We have
\begin{align}
[h_a,K^{0i}{}_i] &= -B_{a0} K^{0i}{}_i\,, & 
(i&=1,2,\ldots,n-1)\,
\end{align}
for all the ${n-1}$ linearly independent elements $K^{0i}{}_i$, and thus they all have
the same weight as $f_0=K^0$. In other words, the
negative $-\alpha_0$ of the simple root $\alpha_0$
has multiplicity $n-1$ in $W(n)$, meaning that the corresponding root space is $(n-1)$-dimensional (whereas $\alpha_0$ itself has
multiplicity one as usual). Furthermore, we note that the adjoint
action of $e_0$ maps the root space of $-\alpha_0$, spanned by all
$K^{0i}{}_i$, injectively to the subspace of the Cartan subalgebra
spanned by $h_0,h_2,h_3\ldots,h_{n-1}$, (that is, the Cartan generators
corresponding to simple roots orthogonal to $\alpha_0$).  Defining
generators
\begin{align}
f_{0i} = ({\rm ad}\,e_0)^{-1} (h_i) = K^{0{(i-1)}}{}_{i-1} - K^{0i}{}_{i}\label{FoiDef}
\end{align}
corresponding to $-\alpha_0$ for $i=2,3,\ldots,n-1$, we thus have $[e_0,f_{0i}]=h_i$. We set $f_{00}
\equiv f_0$ to make this relation valid also for $i=0$. Henceforth,
whenever $f_{0a}$ appears we assume $a=0,2,3,\ldots,r$ (with $r=n-1$ in the case considered here),
and whenever $f_{a}$ appears we assume $a \neq 0$. We thus
have the $\mathbb{Z}_2$-graded set of generators
\begin{align} \label{setofgenerators}
\scr S=\{e_a,f_{0a},f_a\}
\end{align}
where $e_0$ and $f_{0a}$ are odd, and the others are even, and we seek
relations generating an ideal $C$ of the free Lie superalgebra $F$ generated by $\scr S$ such that $F/C$
is isomorphic to $W(n)$.

Of course, one way of finding relations that generate the ideal $C$ would be to write all the commutation relations (\ref{gencommrel})
with the basis elements expressed in terms of the generators $e_a,f_{0a},f_a$. But there would be a great deal of redundancy in
such a set of relations, and, most importantly, they would be applicable only to the case of $A(n-1,0)$,
not to the other cases we are interested in. Thus we seek relations that are more fundamental in this sense.
The obvious starting point is the Chevalley--Serre relations for $\scr B$ that do not involve $f_{0}$,
and the relations
\begin{align}
[e_0,f_{0a}]&=h_a\,,
& 
[h_a,f_{0b}]&=-B_{a0}f_{0b}
\end{align}
explained above. In addition we choose the relations
(\ref{e1f0a})--(\ref{eifjf0a}) below,
which can easily be verified in the $A(n-1,0)$ case, as fundamental ones. We will justify this choice for a general $\fg$ in the next section,
and present the final result for $\fg=A_{n-1}$ in Section \ref{IdealJSection}.

\section{The Lie superalgebras
\texorpdfstring{$W(\mathfrak{g})$}{W(g)} and \texorpdfstring{$S(\mathfrak{g})$}{S(g)}}\label{WgSection}
\subsection{Definition of \texorpdfstring{$\widetilde W(\mathfrak{g})$}{W(g)}}

We return to the general case where a grey node (node 0) is connected to node 1 in the Dynkin diagram
of a simple and simply laced Kac--Moody algebra $\fg$ of rank $r$, and consider the corresponding $\mathbb{Z}_2$-graded set
of generators $\scr S=\{e_a,f_{0a},f_a\}$.
Motivated by the observations made for $W(n)$ at the end of the preceding section, we
define $\widetilde W=\widetilde W(\mathfrak{g})$ as the Lie superalgebra generated
by the set $\scr S$ 
modulo the relations
\begin{align}
[h_a,e_b]=B_{ab}e_b\,,  \qquad   
[h_a,f_b]=-B_{ab}f_b\,, \qquad    
[e_a,f_b]=\delta_{ab}h_b\,,  \label{eigen}
\end{align}
\begin{align}
({\rm ad}\,e_a)^{1-B_{ab}}(e_b)=({\rm ad}\,f_a)^{1-B_{ab}}(f_b)=0\,,
\label{serre}
\end{align} 
\begin{align}
[e_0,f_{0a}]&=h_a\,, \label{e0f0a}
\end{align}
\begin{align} \label{haf0b}
[h_a,f_{0b}]&=-B_{a0}f_{0b}\,,
\end{align}
\begin{align}
[e_1,f_{0a}]&=0\,, \label{e1f0a}
\end{align}
\begin{align}
[e_a,[e_a,f_{0b}]]=[f_a,[f_a,f_{0b}]]=0\,, \label{dubbel}
\end{align}
\begin{align}
[e_i,[f_j,f_{0a}]]&= \delta_{ij}B_{aj}f_{0j}\,, \label{eifjf0a}
\end{align}
for $i,j=2,3,\ldots,r$, where we have defined $h_a$ by (\ref{e0f0a}) for $a\neq 1$, and $h_1=[e_1,f_1]$
(recall that $a \neq 1$ in $f_{0a}$, and $a\neq 0$ in $f_a$).

There is some redundancy in (\ref{eigen})--(\ref{eifjf0a}): 
If $a,b=2,3,\ldots,r$, then the last relation in
(\ref{eigen}) follows by acting with $e_0$ on (\ref{eifjf0a}), and
acting with $e_1$ on (\ref{e1f0a}) gives the relation $[e_a,[e_a,f_{0b}]]=0$ in
(\ref{dubbel}) for $a=1$. Furthermore, for $a=0$ this relation follows from
$[e_0,[e_0,e_1]]=0$ in (\ref{serre}) since
\begin{align}
4[e_0,[e_0,f_{0b}]]&=2[[e_0,e_0],f_{0b}]=[[f_1,[e_1,[e_0,e_0]]],f_{0b}]\nn\\
&=-2[[f_1,[e_0,[e_0,e_1]]],f_{0b}]\,.
\end{align}
The different relations (\ref{dubbel}) for different values of $b$ are not independent either. For example, if $B_{23}=-1$, then we have
\begin{align}
[f_2,[f_2,f_{03}]]&=-[f_2,[f_2,[e_3,[f_3,f_{02}]]]]\nn\\
&=-[e_3,[f_2,[f_2,[f_3,f_{02}]]]]\nn\\
&=-2[e_3,[f_2,[f_3,[f_2,f_{02}]]]]+[e_2,[f_3,[f_2,[f_2,f_{02}]]]]\nn\\
&=-2[f_2,[h_3,[f_2,f_{02}]]]+[e_2,[f_3,[f_2,[f_2,f_{02}]]]]\nn\\
&=-2[f_2,[f_2,f_{02}]]+[e_2,[f_3,[f_2,[f_2,f_{02}]]]]\,,
\end{align}
where we have used $[f_2,[f_2,f_3]]=0$, expanded as in (\ref{serre-expanded}).
Thus (\ref{dubbel}) could be replaced by
\begin{align}
[e_i,[e_i,f_{02}]]&=[f_a,[f_a,f_{02}]]=0\,,\nn\\
[e_i,[e_i,f_{00}]]&=[f_a,[f_a,f_{00}]]=0
\end{align}
for $i=2,3,\ldots,r$. Also in (\ref{e0f0a})--(\ref{e1f0a}) it is sufficient to consider $f_{02}$ and $f_{00}$.

Note the absence of 
relations setting
$[e_i,f_{0a}]$ and $[f_i,f_{0a}]$ to zero for any
$i=2,3,\ldots,r$. However, if $B_{ai}=0$, then it follows from
(\ref{dubbel}) and (\ref{eifjf0a}) that $[e_i,f_{0a}]=[f_i,f_{0a}]=0$.

We will proceed under the assumption that the algebra $\widetilde W(\fg)$ is non-trivial, \ie, that the relations (\ref{eigen})--(\ref{eifjf0a})
do not generate the whole free Lie superalgebra $F$ generated by $\scr S$. We have no proof that this assumption is true for {arbitrary} $\fg$,
only for finite-dimensional $\fg$ (as we will see later in this section), and for $\fg=E_r$ with $r\geq 11$ and the ``outermost'' node as node 1
(as we will show in Section \ref{WEn-section}). In the $A$ case, we have already seen that the non-trivial Lie superalgebra
$W(n)$ satisfies the relations. 

\subsection{The \texorpdfstring{$\mathbb{Z}$}{Z}-grading on \texorpdfstring{$\widetilde W(\mathfrak{g})$}{W(g)} and
definition of \texorpdfstring{$W(\mathfrak{g})$}{W(g)}}
The free Lie superalgebra $F$ generated by $\scr S$
is spanned by all elements
\begin{align} \label{F-basis}
[x_1,[x_2,\ldots,[x_{p-1},x_p]\cdots]] \qquad (p\geq 1),
\end{align}
where each $x_i$ ($i=1,2,\ldots,p$) is equal to $e_a$, $f_a$ or $f_{0a}$. It follows that $F$ has a consistent $\mathbb{Z}$-grading where $F_k$
is spanned by all elements of the form (\ref{F-basis}) such that, among the generators $x_1,\ldots,x_p$,
the number of $e_0$ minus the number of $f_{0a}$
is equal to $k$. Thus, if $F_{(i,j)}$ is the subspace of $F$ spanned by all elements of the form (\ref{F-basis}) where
$e_0$ appears $i$ times and $f_{0a}$ appears $j$ times (in total, for possibly different $a$), then
\begin{align}
F_k = \bigoplus F_{(i,j)}\,, \label{Z-gradingonF}
\end{align}
where the sum ranges over all pairs $(i,j)$ of non-negative integers such that $i-j=k$.

Let $I$ be the ideal generated by the relations (\ref{eigen})--(\ref{eifjf0a}). 
Since these relations are homogeneous with respect to the $\mathbb{Z}$-grading (\ref{Z-gradingonF}) on $F$, this $\mathbb{Z}$-grading is preserved
in the quotient $\widetilde W =F/I$, and we have
\begin{align}
\widetilde W_k = \sum \widetilde W_{(i,j)}\,, \label{Z-gradingonWtilde}
\end{align}
where, as before,  the sum ranges over all pairs $(i,j)$ of non-negative integers such that $i-j=k$.

Note that the sum in (\ref{Z-gradingonWtilde}) is not (necessarily) direct, since the Lie superalgebra $\widetilde W$ is not free. In fact, we will show that
$\widetilde W_{(i,j)}=\widetilde W_{(i',j')}$ for all pairs $(i,j)$ and $(i',j')$ such that $i-j=i'-j'$,
or equivalently
\begin{align}
[\widetilde W_{(1,0)}, \widetilde W_{(0,1)}] = \widetilde W_{(0,0)}\,. \label{Z-gradingcondition}
\end{align}
In particular $\widetilde W_k = \widetilde W_{(k,0)}$
and $\widetilde W_{-k} = \widetilde W_{(0,k)}$ for all $k \geq 0$.

Since (\ref{Z-gradingonWtilde}) is a $\mathbb{Z}$-grading,
$\widetilde W_{-1} \oplus \widetilde W_{0} \oplus \widetilde W_1$ is a local Lie superalgebra. We define $W(\mathfrak{g})$ as the minimal Lie superalgebra with
this local part. Thus we can write, as usual, $\widetilde W_i = W_i$ for $i=0,\pm1$.
Note however that $\widetilde W(\fg)$ is not the maximal Lie superalgebra with local part $\widetilde W_{-1} \oplus \widetilde W_{0} \oplus \widetilde W_1$,
since we, for convenience, have included the relation $[e_0,[e_0,e_1]]=0$ in (\ref{serre}), which occurs at level 2.

We define $\widetilde S=\widetilde S(\fg)$ as the subalgebra of $\widetilde W(\fg)$ generated by the subset $\scr S\backslash \{f_{00}\}$,
and $S=S(\fg)$
as the minimal Lie superalgebra with local part $\widetilde S_{-1} \oplus 
\widetilde S_{0} \oplus \widetilde S_{1}$ with respect to the $\mathbb{Z}$-grading inherited from $\widetilde W(\fg)$.
From our results for $W(\fg)$ below, it is straightforward to deduce the corresponding results for $S(\fg)$. We will not do that explicitly, but focus
on $W(\fg)$. We will come back to $S(\fg)$ later when we relate the results to the construction of the tensor hierarchy algebra in \cite{Palmkvist:2013vya}. 

\subsection{The subalgebra \texorpdfstring{$\widetilde W'$}{W'}}
Obviously, the generators associated to node 1 play a distinguished role in the relations 
(\ref{eigen})--(\ref{eifjf0a}). Therefore, it is convenient to first study the subalgebra of $\widetilde W$
generated by 
the subset $\scr S'=\scr S\backslash \{e_1,f_1\}$ of $\scr S$.
We denote this subalgebra of $\widetilde W$ by $\widetilde W'$ and consider the $\mathbb{Z}$-grading inherited from $\widetilde W$,
which in turn originates from the $(\mathbb{N} \times \mathbb{N})$-grading of the free Lie superalgebra $F$.
Thus $\widetilde {W'}_{(0,0)}$ is the subalgebra of $\widetilde W$ generated by $e_i,f_i,h_0$ for $i=2,3,\ldots,r$,
and since these generators commute with $e_0$, the subspace $\widetilde {W'}_{{(1,0)}}$ is one-dimensional (spanned by $e_0$). Since furthermore 
$[e_0,f_{0a}]=h_a \in \widetilde {W'}_{(0,0)}$,
it follows that the condition for $\widetilde W'$ corresponding to (\ref{Z-gradingcondition}) holds,
\begin{align}
[\widetilde {W'}_{(1,0)}, \widetilde {W'}_{(0,1)}] = \widetilde {W'}_{(0,0)}\,. \label{Z-gradingcondition'}
\end{align}
Thus $\widetilde W'{}_k = \widetilde W'{}_{(k,0)}$
and $\widetilde W'{}_{-k} = \widetilde W'{}_{(0,k)}$ for all $k \geq 0$, in particular for $k=0,1$.

Having described $\widetilde W'{}_1$ as the one-dimensional subspace of $\widetilde W'$ spanned by $e_0$, we now proceed to $\widetilde W'{}_0$,
and in the next subsection to $\widetilde W'{}_{-1}$.
\begin{proposition}The subalgebra $\widetilde W'{}_0=\widetilde W'{}_{(0,0)}$ of $\widetilde W$
is isomorphic to $\fg' \oplus \mathbb{K}$.
\end{proposition}
\Pf It is clear that $\fg' \oplus \mathbb{K}$ can be constructed as the Lie superalgebra generated by $\scr S'=\{e_0,f_{0a},e_i,f_i\}$
for $i=2,3,\ldots,r$ modulo the relations among (\ref{eigen})--(\ref{eifjf0a}) that only involve these generators.
Thus there is a surjective
homomorphism from $\fg'$ to $\widetilde {W'}_{0}$. The kernel of this homomorphism must be an ideal of $\fg'\oplus \mathbb{K}$. Since we assume
that $\fg'$ is simple, the ideal is either $\fg'$ or the one-dimensional center or zero. It is easy to see that in the first two cases
the whole of $\widetilde W$ would be zero, contradicting the assumption above.
We conclude that the homomorphism is injective, and thus an isomorphism.
\qed

Consider the map $\scr S' \to \scr S'$ given by
\begin{align} \label{W'auto}
e_i &\mapsto \pm f_i\,, & f_i &\mapsto  \pm e_i\,, & e_0 &\mapsto e_0\,, & f_{0a} &\mapsto -f_{0a}\,.
\end{align}
Since those of the 
relations (\ref{eigen})--(\ref{eifjf0a}) that do not involve $e_1$ or $f_1$ are invariant under this map,
it can be extended to an automorphism of the Lie superalgebra $\widetilde W'$.

We will now state and prove some additional relations that hold in $\widetilde W(\fg)'$, as consequences of the defining ones (\ref{eigen})--(\ref{eifjf0a}).
Given the automorphism (\ref{W'auto}), these relations always come in pairs, and we will only prove half of them explicitly. The other half then
follow by applying the automorphism. 

\begin{proposition}\label{additionalrelations}
For $i,j=2,3,\ldots,r$ and $i\neq j$, the following relations hold in $\widetilde W'$:
\begin{align} \label{serre-like}
({\rm ad}\,e_i)^{1-B_{ij}}({\rm ad}\,e_j)(f_{0a})=({\rm
    ad}\,f_i)^{1-B_{ij}}({\rm ad}\,f_j)(f_{0a})=0\,, 
\end{align}
\begin{align} \label{proportionality}
B_{ia}[e_i,f_{0b}]&=B_{ib}[e_i,f_{0a}]\,,&B_{ia}[f_i,f_{0b}]&=B_{ib}[f_i,f_{0a}]\,.
\end{align}
\end{proposition}
\Pf
If $B_{ij}=0$, we have
\begin{align}
  0&=\text{ad}\,[e_i,e_j]=[\text{ad}\,e_i,\,\text{ad}\,e_j]
  =\text{ad}\,e_i\,\text{ad}\,e_j-\text{ad}\,e_j\,\text{ad}\,e_i\,, 
\end{align}
and we get
\begin{align}
2[e_i,[e_j,f_{0a}]]&=2[e_i,[e_j,f_{0a}]]+[f_i,[e_j,[e_i,[e_i,f_{0a}]]]]\nn\\
&=[h_i,[e_i,[e_j,f_{0a}]]]+[f_i,[e_i,[e_i,[e_j,f_{0a}]]]]\nn\\
&=[e_i,[f_i,[e_i,[e_j,f_{0a}]]]]\nn\\
&=[e_i,[h_i,[e_j,f_{0a}]]]+[e_i,[e_i,[f_i,[e_j,f_{0a}]]]]=0\,.
\end{align}
Likewise, if $B_{ij}=-1$, we have
\begin{align} \label{serre-expanded}
  0&=\text{ad}\,[e_i,[e_i,e_j]]
  =[\text{ad\,}e_i,\text{ad}\,[e_i,e_j]]=[\text{ad}\,e_i,[\text{ad}\,e_i,
      \text{ad}\,e_j]]\nn\\
  &=\text{ad}\,e_i\,\text{ad}\,e_i\,\text{ad}\,e_j-2\,\text{ad}\,e_i\,
  \text{ad}\,e_j\,\text{ad}\,e_i\,
+\text{ad}\,e_j\,\text{ad}\,e_i\,\text{ad}\,e_i\,,
\end{align}
and we get
\begin{align} \label{serre-like-proof}
  3[e_i,[e_i,[e_j,f_{0a}]]]&=3[e_i,[e_i,[e_j,f_{0a}]]]
  +4[f_i,[e_i,[e_j,[e_i,[e_i,f_{0a}]]]]]\nn\\
  &\quad\,-2[f_i,[e_j,[e_i,[e_i,[e_i,f_{0a}]]]]]
  -[f_i,[e_i,[e_j,[e_i,[e_i,f_{0a}]]]]]\nn\\
&=[h_i,[e_i,[e_i,[e_j,f_{0a}]]]]+2[f_i,[e_i,[e_i,[e_j,[e_i,f_{0a}]]]]]\nn\\
&\quad\,-[f_i,[e_i,[e_j,[e_i,[e_i,f_{0a}]]]]]\nn\\
&=[h_i,[e_i,[e_i,[e_j,f_{0a}]]]]+[f_i,[e_i,[e_i,[e_i,[e_j,f_{0a}]]]]]\nn\\
&=[e_i,[f_i,[e_i,[e_i,[e_j,f_{0a}]]]]]\nn\\
&=-\,[e_i,[h_i,[e_i,[e_j,f_{0a}]]]]+[e_i,[e_i,[f_i,[e_i,[e_j,f_{0a}]]]]]\nn\\
&=-\,[e_i,[h_i,[e_i,[e_j,f_{0a}]]]]-[e_i,[e_i,[h_i,[e_j,f_{0a}]]]]]\nn\\
&\quad\,+[e_i,[e_i,[e_i,[f_i,[e_j,f_{0a}]]]]]=0\,.
\end{align}
Finally we have
\begin{align} \label{prop2}
B_{ia}[f_i,f_{0i}]&=[f_i,[e_i,[f_i,f_{0a}]]]\nn\\
&=-\,[h_i,[f_i,f_{0a}]]-[e_i,[f_i,[f_i,f_{0a}]]]=2[f_i,f_{0a}]\,,
\end{align}
where $a$ can be replaced by $b$. 
Thus (\ref{proportionality})
follows directly if $B_{ia}$ or $B_{ib}$ is zero, and otherwise by 
\begin{align}
[f_i,f_{0i}]=\frac{2}{B_{ia}}[f_i,f_{0a}]=\frac{2}{B_{ib}}[f_i,f_{0b}]\,.
\end{align}
\qed

The meaning of this proposition is that we lose the $r$-fold multiplicity of the roots at level $-1$ when we apply 
 ${\rm ad}\,e_i$ or ${\rm ad}\,f_i$ to $f_{0a}$, for $i=2,3,\ldots,r$. As we have seen, $-\alpha_0$ has multiplicity $r$, whereas the proposition says that
the roots $-\alpha_0+\alpha_i$
have multiplicity one, since for different $a$, all root vectors $[e_i,f_{0a}]$ are proportional to each other
(and correspondingly for $-\alpha_0-\alpha_i$ and $[f_i,f_{0a}]$).

\subsection{Determining \texorpdfstring{$\widetilde {W'}_{-1}$}{W'-1} when \texorpdfstring{$\fg'$}{g'} is finite-dimensional}
From now on, we assume that $\fg'$ is finite-dimensional (and still simply laced). Then we know that all roots
$\alpha$ of $\fg'$ have length squared $(\alpha,\alpha)=2$ with our normalization, and that the sum $\alpha+\beta$
of two positive roots $\alpha,\beta$ is a root if and only if $(\alpha,\beta)=-1$. We will use this repeatedly in the following proofs.
At the end of this section we will comment on the case when $\fg'$ is infinite-dimensional.

For all positive roots $\alpha$ of $\fg'$, and $i=2,3,\ldots,r$, let 
$\{e_\alpha,f_\alpha,h_i\}$
be a Chevalley basis of $\fg$, so that
\begin{align}
[e_\alpha,f_\alpha]&=h_\alpha\,, &
[h_\alpha,e_\beta]&=(\alpha,\beta)e_\beta\,, &
[h_\alpha,f_\beta]&=-(\alpha,\beta)f_\beta\,,
\end{align}
where $h_\alpha=\sum_{i=2}^r a_ih_{i}$ if
$\alpha=\sum_{i=2}^ra_i\alpha_i$.

We introduce the following notation: for $\alpha=\sum_{i=2}^ra_i\alpha_i$, set $f_{0\alpha}=\sum_{i=2}^r a_if_{0i}$. Then 
$f_{0\alpha_i}=f_{0i}$, and the relation (\ref{eifjf0a}) can be written
\begin{align}
[e_i,[f_j,f_{0\alpha}]]&= \delta_{ij}(\alpha,\alpha_j)f_{0j}\,. \label{eifjf0a2}
\end{align}

\begin{lemma} \label{tredelat}
For any positive root $\alpha$ of $\mathfrak{g}'$ and for $i=2,3,\ldots,r$, we have
\begin{align}
(\alpha_i,\alpha)\leq 0 &\qquad\Rightarrow\qquad [f_i,[e_\alpha,f_{0a}]]=0\,, 
\label{lemma-a}\\
(\alpha_i,\alpha)\geq0 &\qquad\Rightarrow\qquad  [e_i,[e_{\alpha},f_{0a}]]=0\,, 
\label{lemma-c}
\end{align}
\begin{align}
[e_{\alpha},[e_i,f_{0a}]]=B_{ia}[e_i,[e_\alpha,f_{0i}]] \label{lemma-b}\,.
\end{align}

\end{lemma}

\Pf
First we note that (\ref{lemma-c}) is a consequence of
(\ref{lemma-b}), since if $(\alpha_i,\alpha)\geq0$, then 
$[e_i,e_\alpha]=0$ and
\begin{align}
[e_{\alpha},[e_i,f_{0a}]]=[e_i,[e_\alpha,f_{0i}]]\,.
\end{align}
Together with (\ref{lemma-b}) this gives (\ref{lemma-c}). Thus it
suffices to prove (\ref{lemma-a}) and (\ref{lemma-b}), 
which can be done by induction
over the height of $\alpha$, denoted ${\rm ht}\,\alpha$. If ${\rm
  ht}\,\alpha =1$, then $\alpha$ 
is a simple root, say $\alpha=\alpha_k$. Then (\ref{lemma-a}) and
(\ref{lemma-b}) follow from  
(\ref{eifjf0a})
and 
\begin{align} \label{basfall-b}
[e_k,[e_i,f_{0a}]]&=\frac{B_{ia}}2[e_k,[e_i,f_{0i}]]=
\frac12[e_k,[e_i,[e_i,[f_i,f_{0a}]]]]\nn\\
&=-\frac12[e_i,[e_i,[e_k,[f_i,f_{0a}]]]]+[e_i,[e_k,[e_i,[f_i,f_{0a}]]]]\nn\\
&=[e_i,[e_k,[e_i,[f_i,f_{0a}]]]]=B_{ia}[e_i,[e_k,f_{0i}]]\,,
\end{align}
respectively.

Suppose now that the lemma holds for roots $\alpha'$ such that ${\rm
  ht}\,\alpha'\leq p$ for some $p\geq 1$.
Thus the induction hypothesis consists of the three parts
\begin{align}
(\alpha_i,\alpha')\leq 0 &\qquad\Rightarrow\qquad [f_i,[e_{\alpha'},f_{0a}]]=0\,, 
\label{lemma-a'}\\
(\alpha_i,\alpha')\geq0 &\qquad\Rightarrow\qquad  [e_i,[e_{\alpha'},f_{0a}]]=0\,, 
\label{lemma-c'}
\end{align}
\begin{align}
[e_{\alpha'},[e_i,f_{0a}]]=B_{ia}[e_i,[e_{\alpha'},f_{0i}]] \label{lemma-b'}\,,
\end{align}
corresponding to (\ref{lemma-a})--(\ref{lemma-b}).
    
  Any root $\alpha$ with ${\rm
  ht}\,\alpha = p+1$ can be written $\alpha=\alpha'+\alpha_j$,
where $(\alpha,\alpha_j)=-1$, the condition $(\alpha_i,\alpha)\leq 0$
implies $\alpha',\alpha_j\neq \alpha_i$, and with
$e_\alpha=[e_j,e_{\alpha'}]$ we get
\begin{align}
  [f_i,[e_\alpha,f_{0a}]]&=[f_i,[e_j,[e_{\alpha'},f_{0a}]]]
  -[f_i,[e_{\alpha'},[e_j,f_{0a}]]]\nn\\
  &=[e_j,[f_i,[e_{\alpha'},f_{0a}]]]-[f_i,[e_{\alpha'},[e_j,f_{0a}]]]\,. 
\end{align}
If now $(\alpha_i,\alpha')\leq 0$, then the first term is zero by the
induction hypothesis, and the second term is equal to 
$[e_{\alpha'},[f_i,[e_j,f_{0a}]]]$,
which is zero by (\ref{eifjf0a}). If $(\alpha_i,\alpha')=1$, then we
can write $\alpha'=\alpha_i+\alpha''$ 
and $e_{\alpha'}=[e_i,e_{\alpha''}]$, where
$(\alpha_i,\alpha'')=-1$. Furthermore, we must have
$(\alpha_i,\alpha_j)=-1$, since $(\alpha_i,\alpha_j)=0$ would imply
\begin{align}
(\alpha_i,\alpha)=(\alpha_i,\alpha_j+\alpha_i+\alpha'')=1\,,
\end{align}
contradicting $(\alpha_i,\alpha)\leq 0$. This means that
\begin{align}
  (\alpha_j,\alpha'')=(\alpha_j,\alpha'-\alpha_i)
  =(\alpha_j,\alpha')-(\alpha_j,\alpha_i)=0\,,
\end{align}
and we get
\begin{align}
  [f_i,[e_\alpha,f_{0a}]]&=[f_i,[e_j,[e_i,[e_{\alpha''},f_{0a}]]]]
  -[f_i,[e_j,[e_{\alpha''},[e_i,f_{0a}]]]]\nn\\
  &\quad\,-[f_i,[e_i,[e_{\alpha''},[e_j,f_{0a}]]]]
  +[f_i,[e_{\alpha''},[e_i,[e_j,f_{0a}]]]]\nn\\
&=-[e_j,[h_i,[e_{\alpha''},f_{0a}]]]+[e_j,[e_i,[f_i,[e_{\alpha''},f_{0a}]]]\nn\\
&\quad\,-[e_j,[e_{\alpha''},[f_i,[e_i,f_{0a}]]]]
+[h_i,[e_{\alpha''},[e_j,f_{0a}]]]\nn\\
&\quad\,-[e_{\alpha''},[h_i,[e_j,f_{0a}]]]\,.
\end{align}
The second term
is zero by (\ref{lemma-a'}), and the others
are proportional to 
$[e_j,[e_{\alpha''},f_{0a}]]$, which is zero by (\ref{lemma-c'}). Thus we have shown that 
(\ref{lemma-a}) holds if ${\rm ht}\,\alpha=p+1$. Now (\ref{lemma-b})
for ${\rm ht}\,\alpha=p+1$ and $(\alpha_i,\alpha)\leq0$ 
can be shown in the same way as (\ref{basfall-b})
with $e_k$ replaced by $e_\alpha$.
If $(\alpha,\alpha_i)=1$, then we can write
$e_\alpha=[e_i,e_{\alpha'}]$, where $(\alpha',\alpha_i)=-1$, and we
get 
\begin{align}
  [e_i,[e_\alpha,f_{0a}]]=[e_\alpha,[e_i,f_{0a}]]
  =[e_i,[e_{\alpha'},[e_i,f_{0a}]]]=\frac12[e_i,[e_i,[e_{\alpha'},f_{0a}]]]\,,
\end{align}
which can be shown to vanish in the same way as
$[e_i,[e_i,[e_j,f_{0a}]]]$ in (\ref{serre-like-proof}), using
(\ref{lemma-a'}).
\qed

We will use (\ref{lemma-a}) in the proof of
Proposition~\ref{chev-basis-prop} below, and also the corresponding
statement 
\begin{align}
(\alpha_i,\alpha)\leq 0 &\qquad\Rightarrow\qquad [e_i,[f_\alpha,f_{0a}]]=0\,,
\end{align}
which follows by applying the automorphism given by (\ref{W'auto}).

\begin{theorem}
The $\mathfrak{g}'$-module $W_{-1}'$ generated by the elements $f_{0a}$ is isomorphic to $\mathfrak{g}'$
itself (\ie, the adjoint module). 
\end{theorem}

\Pf
Consider the linear map $\varphi:\mathfrak{g}' \to W_{-1}'$ given by
\begin{align}
&\varphi: & e_\alpha &\mapsto -\frac1{{\rm ht}\,\alpha}[e_\alpha,f_{0\rho}]\,, &
f_\alpha &\mapsto \frac1{{\rm ht}\,\alpha}[f_\alpha,f_{0\rho}]\,, &
h_i &\mapsto f_{0i}\,, \label{varphidef}
\end{align}
where $\rho$ is the Weyl vector of $\mathfrak{g}'$, satisfying
$(\rho,\alpha_i)=1$ for any $i=2,3,\ldots,r$. 
Then $\varphi$ is injective, with inverse $\ad e_0$.
Furthermore,
\begin{align}
\varphi([e_i,e_\alpha]) &= -\frac1{{\rm
    ht}\,\alpha+1}\bigg([e_i,[e_\alpha,f_{0\rho}]] 
-[e_\alpha,[e_i,f_{0\rho}]]\bigg)\nn\\
&= -\frac1{{\rm ht}\,\alpha+1}\bigg([e_i,[e_\alpha,f_{0\rho}]]
-[e_i,[e_\alpha,f_{0i}]]\bigg)\nn\\
&= -\frac1{{\rm ht}\,\alpha+1}[e_i,[e_\alpha,f_{0(\rho-\alpha_i)}]]
\nn\\
&= -\frac1{{\rm ht}\,\alpha+1}
\frac{(\alpha,\rho-\alpha_i)}{(\alpha,\rho)}[e_i,[e_\alpha,f_{0\rho}]]\,,
\end{align}
which is equal to
\begin{align} \label{intertwiner-lhs}
-\frac1{{\rm ht}\,\alpha+1}\bigg( 1+\frac1{{\rm ht}\,\alpha}
\bigg)[e_i,[e_\alpha,f_{0\rho}]] 
&= -\frac1{{\rm ht}\,\alpha}[e_i,[e_\alpha,f_{0\rho}]]=[e_i,\varphi(e_\alpha)]
\end{align}
if $(\alpha,\alpha_i)=-1$. If $(\alpha,\alpha_i)\geq0$, then
(\ref{intertwiner-lhs}) vanishes by (\ref{lemma-c}), and is thus equal
to $\varphi([e_i,e_\alpha])=0$.  Similarly, we get
$\varphi([f_i,e_\alpha])=[f_i,\varphi(e_\alpha)]$ (both if $\alpha$ is
a simple root and otherwise), $\varphi([e_i,h_j])=[e_i,\varphi(h_j)]$
and $\varphi([f_i,h_j])=[f_i,\varphi(h_j)]$. Thus $\varphi$ is an
isomorphism of $\mathfrak{g}'$-modules, and by repeated use of the
homomorphism property $\varphi([x,y])=[x,\varphi(y)]$, where $x$ is
equal to $e_i$ or $f_i$, it follows that $\varphi$ is surjective. We
conclude that $\varphi$ is an isomorphism.  \qed

\begin{proposition} \label{chev-basis-prop}
For any positive root $\alpha$ of $\fg'$, we have
\begin{align}
[e_\alpha,[f_{\alpha},f_{0a}]]=(\alpha,\alpha_a)f_{0\alpha}\,.
\end{align}
\end{proposition}
\Pf As in Lemma \ref{tredelat}, this can be proven by
induction on the height of $\alpha$. If $\alpha$ is a simple root,
the proposition follows from (\ref{eifjf0a}).  Suppose now that the
proposition holds for roots $\alpha'$ such that ${\rm ht}\,\alpha'\leq
p$ for some $p\geq 1$. Any root $\alpha$ with ${\rm ht}\,\alpha = p+1$
can then be written $\alpha=\alpha'+\alpha_i$, where
$(\alpha',\alpha_i)=-1$.  We set $e_{\alpha}=[e_i,e_{\alpha'}]$ and
$f_{\alpha}=[f_{\alpha'},f_i]$. Then
\begin{align}
[e_\alpha,[f_{\alpha},f_{0\beta}]]&=[e_i,[e_{\alpha'},[f_{\alpha'},[f_i,f_{0\beta}]]]]
-[e_{\alpha'},[e_i,[f_{\alpha'},[f_i,f_{0\beta}]]]]\nn\\
&\quad\,-[e_i,[e_{\alpha'},[f_i,[f_{\alpha'},f_{0\beta}]]]]
+[e_{\alpha'},[e_i,[f_i,[f_{\alpha'},f_{0\beta}]]]]\nn\\
&=[e_i,[h_{\alpha'},[f_i,f_{0\beta}]]]+[e_i,[f_{\alpha'},[e_{\alpha'},[f_i,f_{0\beta}]]]]\nn\\
&\quad\,-[e_{\alpha'},[f_{\alpha'},[e_i,[f_i,f_{0\beta}]]]]
-[e_i,[f_i,[e_{\alpha'},[f_{\alpha'},f_{0\beta}]]]]\nn\\
&\quad\,+[e_{\alpha'},[h_i,[f_{\alpha'},f_{0\beta}]]]
+[e_{\alpha'},[f_i,[e_i,[f_{\alpha'},f_{0\beta}]]]]\nn\\
&=-\,(\alpha',\alpha_i)[e_i,[f_i,f_{0\beta}]]
-(\alpha_i,\beta)[e_{\alpha'},[f_{\alpha'},f_{0i}]]\nn\\
&\quad\,-(\alpha',\beta)[e_i,[f_i,f_{0\alpha'}]]
-(\alpha_i,\alpha')[e_{\alpha'},[f_{\alpha'},f_{0\beta}]]\nn\\
&=(\alpha_i,\beta)f_{0i}+(\alpha_i,\gamma)f_{0\alpha'}\nn\\
&\quad\,+(\alpha',\beta)f_{0i}+(\alpha',\beta)f_{0\alpha'}\nn\\
&=(\alpha,\beta)(f_{0i}+f_{0\alpha'})=(\alpha,\beta)f_{0\alpha}\,.
\end{align}
\qed

As a consequence of Proposition~\ref{chev-basis-prop} we can
replace not
only $e_k$ by $e_\alpha$ in (\ref{basfall-b}), but also 
$e_i$ by $e_\beta$. This gives the following proposition.
\begin{proposition}\label{tredelat2}
For any positive roots $\alpha,\beta,\gamma$ of $\fg'$ we have
\begin{align}
(\alpha,\beta)\geq0 &\qquad\Rightarrow\qquad  [e_\beta,[e_{\alpha},f_{0\gamma}]]=0\,,
\label{lemma-c2}
\end{align}
and
\begin{align}
[e_\alpha,[e_\beta,f_{0\gamma}]]=(\beta,\gamma)[e_\beta,[e_\alpha,f_{0\beta}]]\,.
\end{align}
\end{proposition}

\subsection{Determining \texorpdfstring{$\widetilde W_{(0,1)}$}{W(0,1)} when \texorpdfstring{$\fg$}{g} is finite-dimensional}
We now step out from $\fg'$ to $\fg$, and from $\widetilde W'$ to $\widetilde W$. Proposition~\ref{tredelat2} can be generalised in the following way.
\begin{proposition}\label{tredelat3}
For any positive roots $\alpha$ of $\fg$ and
$\beta,\gamma$ of $\fg'$ we have
\begin{align}
(\alpha,\beta)\geq0 &\qquad\Rightarrow\qquad  [e_\beta,[e_{\alpha},f_{0\gamma}]]=0\,,
\label{lemma-c3}
\end{align}
and
\begin{align} \label{genprop}
[e_\alpha,[e_\beta,f_{0\gamma}]]=(\beta,\gamma)[e_\beta,[e_\alpha,f_{0\beta}]]\,.
\end{align}
\end{proposition}

\begin{figure}[h]
\begin{center}
{\includegraphics[width=2.8in]{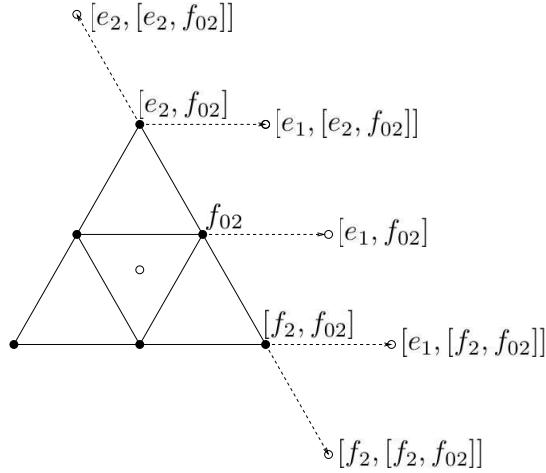}}
\caption{\label{a2null} \it The level $-1$ roots of $S(A_2)$, forming
  the representation $\bf{\bar6}$ of $A_2$. The null
  states of the adjoint (the right edge of the triangle)
  of $A_1$ are shown. All other arrows out of the
level $-1$ root space are obtained via Weyl reflections in the Weyl
group of $A_2$.}
\end{center}
\end{figure}

In particular we can set $e_\alpha=e_1$, which gives $[e_1,[e_\beta,f_{0\gamma}]]=0$, since $[e_1,f_{0\gamma}]=0$.
Thus the annihilation of $f_{0i}=\varphi(h_i)$ upon the
adjoint action of $e_1$ 
(which can be thought of as an arrow in the weight space of $\widetilde W$ pointing out
of the hyperplane which is the weight space of $\mathfrak{g}'$) 
can be transported along any root of $\mathfrak{g}'$, as shown for $\fg=A_2$ in Figure \ref{a2null}.

This can be generalised from $\alpha_1$ to all roots $\alpha_1+\beta$, where $\beta$ is a root of $\fg'$.

\begin{proposition} \label{e1onf0}
Let $\alpha$ be a positive root of $\fg$ such that $(\alpha_0,\alpha)=-1$. Then \linebreak $[e_\alpha,f_{0\beta}]=0$.
\end{proposition}

\Pf We note that the hypothesis means that $\alpha$ appears at level 1 in a level decomposition of $\fg$
with respect to $\alpha_1$. 

Any such root vector $e_\alpha$ can be written as $[e_{\beta_1},[e_{\beta_2},\ldots,[e_{\beta_p},e_1]\cdots]]$ where $\beta_1,\ldots,\beta_p$ are roots
of $\fg'$. Then the proposition can be proven by induction, using Proposition~\ref{tredelat3}.
\phantom{hej}
\qed

We would now like to show that the representation of $\fg$ on $W_{-1}$
is the direct sum of a module dual to $W_1$ and
the irreducible representation of $\fg$ with highest weight
$\Lambda_1+\theta$, where $\theta$ is the highest root of $\fg'$. 
The only further piece of information used in the proof is the
invariance of the relations (\ref{eigen})--(\ref{eifjf0a}) under the
Weyl group of $\fg$, which we now proceed to prove.

\begin{lemma}
The relations (\ref{eigen})--(\ref{eifjf0a}) are preserved by the Weyl
group of $\fg$, with the fundamental Weyl reflections $w_i$ ($i=1,2,\ldots,r$) 
mapping the generators $e_a,f_a,f_{0a},h_a$ to their primed counterparts
\begin{align}
e_a' &=
\begin{cases}
-f_a & \textrm{\rm if }\   i=a \\
[e_i,e_a] & \textrm{\rm if }\   B_{ia}=-1 \\
e_a & \textrm{\rm if }\   B_{ia}=0
\end{cases}\,,
&
f_a' &=
\begin{cases}
-e_a & \textrm{\rm if }\   i=a \\
-[f_i,f_a] & \textrm{\rm if }\   B_{ia}=-1 \\
f_a & \textrm{\rm if }\   B_{ia}=0
\end{cases}\,,\nn\\
f_{0a}' &=
\begin{cases}
f_{0a}-B_{ai}f_{0i} & \textrm{\rm if }\   i\neq1 \\
[f_1,f_{0a}] & \textrm{\rm if }\   i=1
\end{cases} \,, &
h_a' &= h_a-B_{ai}h_i\,. \label{weylgroup}
\end{align} 
\end{lemma}
\Pf
The invariance under the Weyl group of $\fgp$ is quite obvious, and
may easily be checked. The only
additional generator is the Weyl reflection $w_1$ in the hyperplane orthogonal
to $\alpha_1$. The proof proceeds by explicit evaluation of the
identities for the transformed generators.
The identities not containing $f_{0a}$ are the same as in the
BKM superalgebra $\scr B$, and are easily checked.
The relation (\ref{e0f0a}) transforms into
\begin{align}
  [e_0',f_{0i}']&=-[[e_1,e_0],[f_1,f_{0i}]]\nn\\
  &=[e_0,[e_1,[f_1,f_{0i}]]]-[e_1,[e_0,[f_1,f_{0i}]]]\nn\\
  &=[e_0,[h_1,f_{0i}]]-[e_1,[f_1,h_i]]\nn\\
  &=h_i-B_{1i}h_1=h_i'\,,
\end{align}
and (\ref{haf0b}) into
\begin{align}
  [h_a',f_{0i}']&=-[(h_a-B_{a1}h_1),[f_1,f_{0i}]]\nn\\
    &=-(-B_{a1}-B_{a0}-B_{a1}(-2+1))[f_1,f_{0i}]\nn\\
    &=-B_{a0}f_{0i}'\,.
\end{align}
For (\ref{e1f0a}), we obtain
$[e_1',f_{0i}']=[f_1,[f_1,f_{0i}]]=0$.
We will not exhibit all the cases of (\ref{dubbel}) and (\ref{eifjf0a}), only give one example of each. The other cases are similar.
One part of the transformed relation (\ref{dubbel}) is
\begin{align}
  [e_2',[e_2',f_{0i}']]&=-[[e_1,e_2],[[e_1,e_2],[f_1,f_{0i}]]]\nn\\
  &=-[e_1,[e_2,[e_1,[e_2,[f_1,f_{0i}]]]]]
  +[e_1,[e_2,[e_2,[e_1,[f_1,f_{0i}]]]]]\nn\\
  &\quad\,+[e_2,[e_1,[e_1,[e_2,[f_1,f_{0i}]]]]]
  -[e_2,[e_1,[e_2,[e_1,[f_1,f_{0i}]]]]]\nn\\
  &=-[e_1,[e_2,[h_1,[e_2,f_{0i}]]]]+[e_1,[e_2,[e_2,[h_1,f_{0i}]]]]\nn\\
  &\quad\,+[e_2,[e_1,[h_1,[e_2,f_{0i}]]]]-[e_2,[e_1,[e_2,[h_1,f_{0i}]]]]\nn\\
  &=[e_1,[e_2,[e_2,f_{0i}]]]-[e_2,[e_1,[e_2,f_{0i}]]]=0\,,
\end{align}
where we have used
$[e_1,[e_2,f_{0i}]]=0$, see (\ref{genprop}). 
An example of (\ref{eifjf0a}) is
\begin{align}
  [e'_2,[f'_2,f'_{0a}]]&=[[e_1,e_2],[[f_1,f_2],[f_1,f_{0a}]]]\nn\\
  &=[e_1,[e_2,[f_1,[f_2,[f_1,f_{0a}]]]]]]
  -[e_2,[e_1,[f_1,[f_2,[f_1,f_{0a}]]]]]]\nn\\
&=[e_1,[f_1,[h_2,[f_1,f_{0i}]]]]+[e_1,[f_1,[f_2,[e_2,[f_1,f_{0a}]]]]]\nn\\
&\quad\,-[e_2,[h_1,[f_2,[f_1,f_{0a}]]]]-[e_2,[f_1,[f_2,[h_1,f_{0a}]]]]\nn\\
&=-[f_1,[e_2,[f_2,f_{0a}]]]=-B_{2a}[f_1,f_{02}]=B_{2a}f'_{02}\,.
\end{align}
Part of the verification of Weyl invariance
relies on the identity $[f_1,[e_2,f_{0a}]]=0$, which may
be derived from $[f_1,[f_1,[e_1,[e_2,f_{0a}]]]]=0$.
\qed

Thus the map (\ref{weylgroup}) gives an automorphism of $\widetilde W$. The invariance under the Weyl group could be used to prove Proposition 
\ref{chev-basis-prop}, Proposition\ \ref{tredelat2} and Proposition\ \ref{tredelat3}
directly from the corresponding relations for the simple roots
since (in this case) all the positive roots are in one single orbit under the Weyl group \cite{Carbone2014}.
(On the other hand, we used a special case of Proposition 3.7 in the proof of the Weyl group invariance.)
The Weyl group invariance also gives additional parts of Proposition\ \ref{tredelat2} and Proposition\ \ref{tredelat3}, corresponding to (and generalising)
the part (\ref{lemma-a}) of Lemma \ref{tredelat}.

\begin{theorem}
The $\fg$-module $\widetilde W_{0,1}$
is the direct sum of a module dual to $W_1$ and
the irreducible representation of $\fg$ with highest weight
$\lambda=\Lambda_1+\theta$, where $\theta$ is the highest root of $\fg'$.
\end{theorem}
\Pf
The module $R(\Lambda_1)$ dual to $W_{1}$ is generated from its highest weight state
$f_{00}$ as in the
BKM superalgebra $\scr B$.

To prove that $R(\Lambda_1+\theta)$ is a submodule of $\widetilde W_{0,1}$,
consider the element $F_\lambda=[e_{\theta},f_{0j}]$, where $j=2,3,\ldots,r$ is such that
$(\theta,\alpha_j)\neq 0$ (we can always find such a simple root $\alpha_j$).
Then $F_\lambda=[e_{\theta},f_{0j}]$ is 
nonzero since $[f_\theta,[e_{\theta},f_{0j}]]=(\theta,\alpha_j)f_{0\theta}\neq0$.
It is furthermore
a highest
weight vector in the adjoint representation of $\fg'$, and carries $\fg$-weight
$\lambda=\Lambda_1+\theta$. Therefore $F_\lambda$  satisfies
$(\ad f_i)^{\lambda_i+1}F_\lambda=0$ for $i=2,\ldots,r$ \cite{Kac}. 
It is also annihilated by $e_1$, {\it i.e.,}
$[e_1,F_\lambda]=0$. Consider the image of this relation
under a Weyl reflection $w_1$ in the hyperplane orthogonal to
$\alpha_1$. We have $\lambda_1=(\lambda,\alpha_1)=1-c_2$, where $c_2$ is the
Coxeter label of root 2 in $\fg'$. The weight $\lambda$ is thus
orthogonal to $\alpha_1$ if and only if $c_2=1$, which is also a
necessary condition for $\fg$ to be
finite-dimensional. This implies that $w_1(F_\lambda)=F_\lambda$.
We then have $w_1([e_1,F_\lambda])=-[f_1,F_\lambda]=0$. This relation
completes the set of null states for the irreducible
representation of $\fg$ with highest weight $\lambda$, so that
$(\ad f_i)^{\lambda_i+1}F_\lambda=0$ for $i=1,\ldots,r$.

Finally, from $f_{0j}$ we can obtain any $f_{0i}$ with $i=2,3,\ldots,r$ such that $A_{ij}\neq0$
by stepping up and down with generators $e_i$ and $f_i$, and by continuing the procedure  
we can reach all $f_{0i}$ with $i=2,3,\ldots,r$. Thus all such $f_{0i}$ belong to the same module $R(\Lambda_1+\theta)$,
whereas $f_{00}$ belongs to $R(\Lambda_1)$, and there cannot be any other
submodules of $W_{(0,1)}$.
\qed

\subsection{Completing the local Lie superalgebra}

We are now ready to show that the condition (\ref{Z-gradingcondition}) indeed holds when $\fg'$ is finite-dimensional.

\begin{proposition}
We have
\begin{align}
[\widetilde W_{(1,0)},\widetilde W_{(0,1)}]=\widetilde W_{(0,0)}\,.
\end{align}
\end{proposition}
\Pf
We write any element 
in $\widetilde W_{(1,0)}$ as a sum of terms
\begin{align} \label{W10-element}
[e_{\beta_1},[e_{\beta_2},\ldots,[e_{\beta_p},e_0]\cdots]]
\end{align}
where each $\beta_i$ ($i=1,2,\ldots,p$) is a root of $\fg$ such that $(\alpha_0,\beta_i)=-1$,
and any element of $\widetilde W_{(0,1)}$ can be written as a sum of terms
\begin{align}
[x_{1},[x_{2},\ldots,[x_{q},f_{0a}]\cdots]]
\end{align}
where each $x_j$ ($j=1,2,\ldots,q$) is equal to $e_k$ or $f_k$ for $k=1,2,\ldots,r$.
For $q=0$ we have
\begin{align} 
[f_{0a},[e_{\beta_1},[e_{\beta_2},\ldots,[e_{\beta_p},e_0]\cdots]]]&=[e_{\beta_1},[e_{\beta_2},\ldots,[e_{\beta_p},[f_{0a},e_0]]\cdots]]\nn\\
&=[e_{\beta_1},[e_{\beta_2},\ldots,[e_{\beta_p},h_a]\cdots]] \in \widetilde W_{(0,0)} \label{basfall}
\end{align}
by Proposition\ \ref{e1onf0}, and for $q=1$ we get 
\begin{align}
[[x_1,f_{0a}],[e_{\beta_1},[e_{\beta_2},\ldots,[e_{\beta_p},e_0]\cdots]]]&=
[x_1,[f_{0a},[e_{\beta_1},\ldots,[e_{\beta_p},e_0]\cdots]]\nn\\
&\quad\,
-[f_{0a},[x_1,[e_{\beta_1},\ldots,[e_{\beta_p},e_0]\cdots]],
\end{align}
where $[x_1,[e_{\beta_1},\ldots,[e_{\beta_p},e_0]\cdots]]$ in the second term on the right hand side
can be rewritten in the form (\ref{W10-element}) since this is an element of $\widetilde W_{(1,0)}$.
We can then as in (\ref{basfall}) show that each of the two terms on the right hand side belongs to $\widetilde W_{(0,0)}$,
and, continuing in the same way, the proposition follows by induction on $q$.
\qed

To summarise, when $\fg$ is finite-dimensional, the local part $W_{-1} \oplus W_0 \oplus W_1$ of
$W(\fg)$ consists of the $\fg$-modules
\begin{align}
W_{-1} &= R(\Lambda_1 +\theta) \oplus R(\Lambda_1)\,, & W_0 &=\fg \oplus \mathbb{K}\,, & W_1&=\overline{R(\Lambda_1)}\,. \label{localpartofW}
\end{align}
Correspondingly, the local part $S_{-1} \oplus S_0 \oplus S_1$ of
$S(\fg)$ consists of the $\fg$-modules
\begin{align}
S_{-1} &= R(\Lambda_1 +\theta) \,, & S_0 &=\fg\, , & S_1&=\overline{R(\Lambda_1)}\,.
\end{align}
Since all the modules are irreducible in this case, the form of the commutator relations $[S_{-1},S_1]=S_0$ 
is uniquely given by the projection of the tensor product $S_{-1} \otimes S_1$ onto the submodule $S_0$.
Thus, to prove that the definition of $S(\fg)$ here agrees with
the definition of the tensor hierarchy algebra associated to $\fg$ given in \cite{Palmkvist:2013vya}, and reviewed in Section \ref{THA-subsection}, it suffices to check that
the $\fg$-modules in the local parts agree. The only cases where they do not agree are those where 
the module at level $-1$ in the tensor hierarchy algebra is not irreducible but contains a singlet in addition to $S_{-1}$. As explained in
\cite{Palmkvist:2013vya},
this happens precisely
when replacing the grey node in the Dynkin diagram with an ordinary white one gives the Dynkin diagram of an affine Lie algebra, for example
$\fg=E_8$ with the default choice of ``node 1'', where $S_{-1}$ is the module ${\bf 3875}$, but the level $-1$ content of the
tensor hierarchy algebra is ${\bf 3875} + {\bf 1}$. In the application to gauged supergravity the additional singlet is important, but it does not fit naturally
into the algebra from the point of view that we adopt here.

We conclude this section with a few words about the case when $\fg$ is infinite-dimensional.
Then it is no longer true that all roots $\alpha$ of $\fg$ have length squared $(\alpha,\alpha)=2$. In particular there might be a root $\beta$
which satisfies $(\alpha_0,\beta)=-1$ but has length squared $(\beta,\beta)\leq 0$, and thus does not belong to the same Weyl orbit as $e_1$.
Then a corresponding root vector $e_\beta$ does not need to commute with all $f_{0a}$, but $[f_{0a},e_\beta]$ could be a
root vector in $\widetilde W_{(0,1)}$
corresponding to a root $-\alpha_0+\beta$. Taking the commutator
with $e_0$ we get
\begin{align}
[e_0,[f_{0a},e_\beta]]=[h_{0a},e_\beta]-[f_{0a},[e_0,e_\beta]]
\end{align}
where the right hand side now contains a second term in
$\widetilde W_{(1,1)}$ which does not necessarily belong to $\fg=\widetilde W_{(0,0)}$ or vanish.
This is in agreement with the results in \cite{Bossard:2017wxl}, where a tensor hierarchy algebra associated to $E_{11}$
was constructed, and shown to contain elements at level zero beyond the original $E_{11}$ algebra. We will come back to this example in Section
\ref{WEn-section}.

\section{The ideal \texorpdfstring{$J$}{J} of
  \texorpdfstring{$\widetilde{W}(A_{n-1})$}{W(A(n-1))}
\label{IdealJSection}}

Now $W$ can be constructed from $\widetilde W$ as the minimal Lie superalgebra with local part (\ref{localpartofW}),
by factoring out the maximal homogeneous ideal
intersecting the local part trivially. Let $J$ be this ideal of $\widetilde W$. Then $J$ is the direct sum of subspaces $J_k = J \cap \widetilde W_{-k}$
for $k \geq 2$. The intersection of $J$ and $\widetilde W_+$ must be empty since $\widetilde W_+ = \scr B_+$, and $\scr B$ is simple. We conjecture that
$J$ is generated by $J_2$, but we have only proven this for $\fg=A_{n-1}$, the proof of which is the goal of this section.
Throughout the section we assume $\fg=A_{n-1}$ and thus $W=W(\fg)=W(n)$.

\subsection{The intersection between \texorpdfstring{$J$}{J} and
  \texorpdfstring{$\widetilde W_{-2}$}{W(-2)}}
Before stating the identities needed to generate the ideal $J$, 
we examine $\widetilde W_-$, which is
the free Lie superalgebra generated by the subspace ${W}_{-1}$ of $W(A_{n-1})$.
The first observation is that the generators $f_{0a}$ anticommute
among themselves in $W(A_{n-1})$.  Therefore, at
least a part of the ideal $J$ is generated by the anticommutators
$[f_{0a},f_{0b}]$.
It is also straightforward to verify, by acting with $e_0$, that
the $[f_{0a},f_{0b}]$ indeed generate an ideal that intersects the local part trivially.
In order to examine if more generators are
needed, we consider the anticommutator of two elements at level $-1$ in
$\widetilde{W}_{-}$. The level $-1$ generators ($\hat K$) of
$S(A_{n-1})$ form the $A_{n-1}$ representation $(010\ldots01)$. In
$W(A_{n-1})$ there are in addition generators ($K'$) in
$(10\ldots0)$.

All the anticommutators under consideration carry the
$A_{n-1}$ weight $(20\ldots0)$. 
In the freely generated algebra, there are generators at level $-2$
from the anticommutators
\begin{align}
  [\hat K,\hat K]
  &:\,(00010\ldots02)\oplus(0010\ldots01)\oplus(020\ldots02)\cr
  &\quad\oplus(1010\ldots010)\oplus(110\ldots01)\oplus(20\ldots0)\,,\cr
  [\hat K,K']
  &:\,(0010\ldots01)\oplus(010\ldots0)\oplus(110\ldots01)\,,\cr
  [K',K']&:\,(20\ldots0)\,.
  \label{hatKhatK}
\end{align}
Concerning $[\hat K,\hat K]$, the relation 
\begin{align}[f_{0i},f_{0j}]&=0\label{constraint1}\end{align}
 provides ${(n-2)(n-1)}/{2}$ elements at $A_{n-1}$ weight
$(20\ldots0)$. If on the other hand we count the multiplicity of this
 weight in the six representations in $[\hat K,\hat K]$, we
 obtain the numbers in Table \ref{multiplicitytable1}.
 \begin{table}
 \begin{center}
\begin{tabular}{| c  | c |  c | c  | }	

\hline
representation & multiplicity of $(20\ldots0)$     \\ 
\hline
$(00010\ldots02)$ &   $0$  \\ 
$(0010\ldots01)$ &   $0$  \\ 
$(020\ldots02)$ &   $\tfrac{(n-2)(n-1)}{2}$  \\ 
$(1010\ldots010)$ &   $\tfrac{(n-3)(n-2)}{2}-1$  \\ 
$(110\ldots01)$ &   $n-2$  \\ 
$(20\ldots0)$ &   $1$  \\ 
\hline

\end{tabular} 
 \end{center}
 \vspace{4pt}
 \caption{\it \label{multiplicitytable1}Multiplicities of the weight
   $(20\ldots0)$ in some $A_{n-1}$ modules.}
 \vspace{-16pt}
\end{table}
 
The total multiplicity of the weight $(20\ldots 0)$ in $[\hat K,\hat K]$ is $(n-2)(n-1)$. Note that $(0010\ldots01)$,
which is the level $-2$ generator that should be outside the ideal, is
not affected by  relation~(\ref{constraint1}), but also that $(00010\ldots02)$, which should be part of
the ideal, remains untouched by  relation~(\ref{constraint1}).

It remains to verify that  the
remaining elements are eliminated by the  relation~(\ref{constraint1}).
Before checking this, we note that the total multiplicity is larger than
(twice) the number of relations.

Any element in $[\hat K,\hat K]$ with $A_{n-1}$ weight $(20\ldots0)$ comes
from the product of two elements in the adjoint representation of
$A_{n-2}$ at opposite $A_{n-2}$ weights. In addition to the relation
$[f_{0i},f_{0j}]=[\varphi(h_{i}),\varphi(h_{j})]$, one may also consider $[\varphi(e_\alpha),\varphi(f_{\alpha})]$,
where $\alpha$ is a positive $A_{n-2}$ root, and the homomorphism $\varphi$ was defined in (\ref{varphidef}).
The number of positive roots in $A_{n-2}$ is
$\tfrac{(n-2)(n-1)}{2}$. Together with $[f_{0i},f_{0j}]$, the anticommutators
$[\varphi(e_\alpha),\varphi(f_{\alpha})]$ give all the $(n-2)(n-1)$ elements at
$A_{n-1}$ weight $(20\ldots0)$. Now, consider
\begin{align}
  (\ad {e_k})(\ad {f_k})([f_{0i},f_{0j}])
  =2B_{(i|k}[f_{0j)},f_{0k}]-2B_{ik}B_{jk}[\varphi({e_k}),\varphi({f_k})]\,.
\end{align}  
This shows that $[\varphi({e_k}),\varphi({f_k})]$ are not needed as
separate generators in the ideal. Continued action with
$(\ad {e_k})(\ad {f_k})$ gives the full set of $[\varphi(e_\alpha),\varphi(f_{\alpha})]$.

In $[\hat K,K']$, the representation $(110\ldots01)$ should be part
of the ideal. It contains the weight $(20\ldots0)$ with multiplicity
$n-2$, and is set to zero by the relations $[f_{00},f_{0i}]=0$.
In $[K',K']$ the relation $[f_{00},f_{00}]=0$ generates the whole
irreducible $A_{n-1}$ representation $(20\ldots0)$), just as
$[e_0,e_0]=0$ generates its dual at level 2.

At this stage, it remains to remove the representation
$(00010\ldots02)$ in $[\hat K,\hat K]$, and also to relate the
$(0010\ldots01)$'s appearing in $[\hat K,\hat K]$ and
$[\hat K,K']$.
 If, in addition to a pair of
generators $f_{00}$ or $f_{0i}$, one $f_1$ is introduced ($f_{00}$ and
$f_{0i}$ are annihilated by $\ad {e_1}$), the result has $A_{n-1}$
weight $(010\ldots0)$. Counting the multiplicities of this weight in
the representations not already eliminated, we get the numbers in
Table \ref{multiplicitytable2}.

\begin{table}
\begin{center}
\begin{tabular}{| c  | c |  c | c  | }	

\hline
representation & multiplicity of $(010\ldots0)$     \\ 
\hline
$(00010\ldots02)$ &   $\tfrac{(n-4)(n-3)}{2}$  \\ 
$(0010\ldots01)$ &   $n-3$  \\ 
$(010\ldots0)$ &   $1$  \\ 
\hline

\end{tabular} 
\end{center}
\vspace{4pt}
 \caption{\it \label{multiplicitytable2}Multiplicities of the weight
   $(010\ldots0)$ in some $A_{n-1}$ modules.}
\vspace{-16pt}
\end{table}

The element $[f_{0i},[f_1,f_{0j}]]$ generates an
ideal if $i,j=3,\ldots,n-1$, so that nodes $i,j$ are disconnected from
node 1, which is easily seen by commuting with $e_0$. Its symmetric
part vanishes modulo $[f_{0i},f_{0j}]$ by the Jacobi identity. The
antisymmetric part gives exactly the ${(n-4)(n-3)}/{2}$ relations
needed to eliminate $(00010\ldots02)$. The $(0010\ldots01)$, which is
not part of the ideal, contains the $n-3$ elements
$[f_{02},[f_1,f_{0j}]]$, $j=3,\ldots,n-1$.

Finally, one can check that the elements $[(f_{02}-f_{00}),[f_1,f_{0i}]]$ for $i=3,\ldots,n-1$ are
annihilated by $\ad {e_0}$, and generate an ideal that intersects the local part trivially. This provides the
relation between the $(0010\ldots01)$'s
appearing in $[\hat K,\hat K]$ and $[\hat K,K']$.

This completes the examination of the relations at level $-2$. 
We summarise the result.
\begin{theorem}
The intersection $J_2$ of the ideal $J$ and the subspace $\widetilde W_{-2}$ of $\widetilde W$ is generated by the relations
\begin{align}
[f_{0a},f_{0b}]&=0\,,\nn\\
[f_{0i},[f_1,f_{0j}]]&=0\,,\quad i,j\geq3\,,\nn\\
[(f_{02}-f_{00}),[f_1,f_{0i}]]&=0\,,\quad i\geq3\,\label{IdealJ2'}.
\end{align}
\end{theorem}

Before
being in a position to state that this is the full set of generators
of the ideal $J$
at negative levels, we need to show that no new ideals appear at lower
levels.
This can  be done in several ways. Below we will use the
relations recursively and verify that
$W(A_{n-1})$ arises by
repeated use of the relations. 
This proof has the potential of generalization to the $D$-
and $E$-series. 

We thus consider $\widetilde W / J_2$.
We will use repeatedly in the proofs of the two lemmas below that
$[K^{ab}{}_b,K^{cd}{}_d]=0$
in $\widetilde W / J_2$ if $b\neq c$ and $a \neq d$. (As elsewhere in the paper, repeated indices should {\it not} be summed over.)

For $p\geq 3$ and indices $a_1,\ldots,a_p,b$ such that either $b=a_p$
or $b\neq a_1,\ldots, a_p$, 
define $\widetilde K^{a_1\cdots a_p}{}_b$ by $\widetilde K^{a_1a_2a_3}{}_b = K^{a_1a_2a_3}{}_b$ for $p=3$, and
recursively by
\begin{align}
  \widetilde K^{a_1\cdots a_p}{}_b
  &= [K^{a_1a_2}{}_{a_2},\widetilde K^{a_2\cdots a_p}{}_b]\,
\end{align}
for $p\geq 4$.
Let $V_{-p+1}$ be the subspace of $(\widetilde W / J_2)_{-p+1}$
spanned by all such $\widetilde K^{a_1\cdots a_p}{}_b$, 
and set 
\begin{align}
V= \bigg(\bigoplus_{k\geq 3} V_{-k} \bigg)\oplus (\widetilde W /J_2)_{-2} \oplus W_{-1} \oplus W_0 \oplus W_1\,.
\end{align}

\begin{lemma} \label{antisymmetryKtilde}
We have
\begin{align}
\widetilde K^{a_1\cdots a_p}{}_b &= \left\{ \begin{array}{lcc}
\widetilde K^{[a_1\cdots a_p]}{}_b&\, {\rm \it if}\quad&a_p \neq b\,,\\
\widetilde K^{[a_1\cdots a_{p-1}]a_p}{}_b&\, {\rm \it if}\quad&a_p = b\,.
\end{array}\right.
\end{align}
\end{lemma}
\Pf
We prove this by induction. By the definition and the induction hypothesis we then already have
\begin{align}
\widetilde K^{a_1\cdots a_p}{}_b &= \left\{ \begin{array}{lcc}
\widetilde K^{a_1a_2[a_3\cdots a_p]}{}_b&\, {\rm if}\quad&a_p \neq b\,,\\
\widetilde K^{a_1a_2[a_3\cdots a_{p-1}]a_p}{}_b&\, {\rm if}\quad&a_p = b\,,
\end{array}\right.
\end{align}
and it remains to show that
\begin{align}
\widetilde K^{a_1a_2a_3a_4\cdots a_p}{}_b&= - \widetilde K^{a_2a_1a_3a_4\cdots a_p}{}_b=- \widetilde K^{a_3a_2a_1a_4\cdots a_p}{}_b\,.
\end{align}
Suppose $\widetilde K^{a_1\cdots a_p}{}_b \neq 0$. By
the antisymmetry in the upper indices of $K^{a_1 a_2}{}_{a_2}$ and the induction hypothesis we 
can assume that all indices $a_1,a_2,\ldots, a_p$ are
distinct. It follows that $b \neq a_1,a_2,\ldots, a_{p-1}$.
Now we have
\begin{align}
\widetilde K^{a_1a_2a_3\cdots a_p}{}_b &=
       [K^{a_1a_2}{}_{a_2},[K^{a_2a_3}{}_{a_3},\widetilde K^{a_3a_4\cdots
             a_p}{}_b]]\nn\\ 
&=[[K^{a_1a_2}{}_{a_2},K^{a_2a_3}{}_{a_3}],\widetilde K^{a_3a_4\cdots a_p}{}_b]\nn\\
&=-[[K^{a_2a_1}{}_{a_1},K^{a_1a_3}{}_{a_3}],\widetilde K^{a_3a_4\cdots a_p}{}_b]\nn\\
&=-[K^{a_2a_1}{}_{a_1},[K^{a_1a_3}{}_{a_3},\widetilde K^{a_3a_4\cdots
             a_p}{}_b]]=-\widetilde K^{a_2a_1a_3a_4\cdots a_p}{}_b\,,
\end{align}
and similarly
\begin{align}
\widetilde K^{a_1a_2a_3\cdots a_p}{}_b &=
       [K^{a_1a_2}{}_{a_2},[K^{a_2a_3}{}_{a_3},[K^{a_3a_4}{}_{a_4},\widetilde K^{a_4\cdots
             a_p}{}_b]]\nn\\ 
&=[K^{a_1a_2}{}_{a_2},[[K^{a_2a_3}{}_{a_3},K^{a_3a_4}{}_{a_4}],\widetilde K^{a_4\cdots a_p}{}_b]]\nn\\
&=-[K^{a_1a_2}{}_{a_2},[[K^{a_3a_2}{}_{a_2},K^{a_2a_4}{}_{a_4}],\widetilde K^{a_4\cdots a_p}{}_b]]\nn\\
&=[K^{a_3a_2}{}_{a_2},[[K^{a_1a_2}{}_{a_2},K^{a_2a_4}{}_{a_4}],\widetilde K^{a_4\cdots a_p}{}_b]]\nn\\
&=-[K^{a_3a_2}{}_{a_2},[[K^{a_2a_1}{}_{a_1},K^{a_1a_4}{}_{a_4}],\widetilde K^{a_4\cdots a_p}{}_b]]\nn\\
&=-[K^{a_3a_2}{}_{a_2},[K^{a_2a_1}{}_{a_1},[K^{a_1a_4}{}_{a_4},\widetilde K^{a_4\cdots a_p}{}_b]]]=- \widetilde K^{a_3a_2a_1a_4\cdots a_p}{}_b\,.
\end{align}

\qed

\begin{lemma}The subspace $V$ of $\widetilde W/J_2$ is closed under the adjoint action of elements at level $-1$ and $0$, that is,
\begin{align}
[K^{cd}{}_d, \widetilde K^{a_1\cdots a_p}{}_b] \in V\, \label{WClosedMinusOne}
\end{align}
and
\begin{align}
[K^{c}{}_d, \widetilde K^{a_1\cdots a_p}{}_b] \in V\,. \label{WClosedZero}
\end{align}
\end{lemma}

\Pf
Thanks to Lemma \ref{antisymmetryKtilde}, we can assume $c\neq a_2,\ldots,a_p$ and $d \neq a_1,\ldots,a_{p-2}$. Then
\begin{align}
[K^{cd}{}_d, \widetilde K^{a_1\cdots a_p}{}_b] =
[K^{a_1a_2}{}_{a_2},[K^{a_2a_3}{}_{a_3},
      \ldots,[K^{a_{p-2}a_{p-1}}{}_{a_{p-1}},[K^{cd}{}_{d},K^{{a_{p-1}}a_p}{}_{b}]]\cdots]]\,, \nn
\end{align}
where $[K^{cd}{}_{d},K^{{a_{p-1}}a_p}{}_{b}]$ can be written as a linear combination of terms $K^{a_{p-1}ef}{}_g$.
Thus $[K^{cd}{}_d, \widetilde K^{a_1\cdots a_p}{}_b]$ is equal to a corresponding linear combination of terms
\begin{align}
K^{a_1 \cdots a_{p-1}ef}{}_g\,,
\end{align}
and we have proven the first part (\ref{WClosedMinusOne}) of the lemma.
For the second part (\ref{WClosedZero})
we have
\begin{align}
[K^{c}{}_d, \widetilde K^{a_1\cdots a_p}{}_b] =
[K^{c}{}_d, [K^{a_1 a_2}{}_{a_2}  \widetilde K^{a_2\cdots a_p}{}_b]]=
[K^{a_1 a_2}{}_{a_2},[K^{c}{}_d,   \widetilde K^{a_2\cdots a_p}{}_b]]\,, 
\end{align}
and the proof follows by induction, using (\ref{WClosedMinusOne}).
\qed

Since the $A_{n-1}$ module $W_{-1}$ is generated by $K^{cd}{}_d$ for any nonzero $K^{cd}{}_d$
it follows that $[W_{-1},\widetilde{V}_{-p+1}]  \in \widetilde{V}_{-p} $, and
thus 
$\widetilde W/J_2=\widetilde V$. On the other hand, there is an isomorphism
$\widetilde V \leftrightarrow W$ given by 
$\widetilde K^{a_1 \cdots a_p}{}_b \leftrightarrow K^{a_1 \cdots
  a_p}{}_b$, and thus $\widetilde W/J_2$ is isomorphic to $W$. 

We have proven the following theorem.

 \begin{theorem} The ideal $J$ of $\widetilde{W}(A_{n-1})$ is generated by the relations 
 \begin{align}
[f_{0a},f_{0b}]&=0\,,\nn\\
[f_{0i},[f_1,f_{0j}]]&=0\,,\quad i,j\geq3\,,\nn\\
[(f_{02}-f_{00}),[f_1,f_{0i}]]&=0\,,\quad i\geq3\,\label{IdealJ2}.
\end{align}
\end{theorem}

 This gives us the main result of the paper:
  \begin{theorem}\label{maintheorem} The Lie superalgebra $\widetilde{W}(A_{n-1})/J$ is 
  isomorphic to $W(n)$. Thus 
    $W(n)$ has
    generators (\ref{setofgenerators}) 
        and defining relations
(\ref{eigen})--(\ref{eifjf0a}) and (\ref{IdealJ2}).

\end{theorem}

\section{Comments on the \texorpdfstring{$D$}{D} and
  \texorpdfstring{$E$}{E} cases\label{DESection}} 

The tensor hierarchy algebras $W(E_r)$ and $S(E_r)$ are relevant in
exceptional geometry and $W(D_r)$ and $S(D_r)$ are their analogs
in double geometry. The Dynkin diagrams of these algebras, \ie, those of the BKM superalgebras $\scr B(D_r)$ and $\scr B(E_r)$ are given
in Figure \ref{WDEdiagrams}.

\begin{figure}
\begin{center}
\ddiagram 
\end{center}
\begin{center}
\ediagram 
\end{center}
\caption{\it \label{WDEdiagrams}Dynkin diagrams for $\scr B(D_r)$ and $\scr B(E_r)$.}
\end{figure}

The definition of $W(\fg)$ given in Section \ref{WgSection} is
formulated entirely in terms of generators, and holds for any $\fg$.
This means that the local superalgebra $W_{-1}\oplus\ W_0\ \oplus\ W_{1}$ of 
$W(\fg)$ is obtained from the definition of $\widetilde W(\fg)$. The algebra $W(\fg)$ is
then defined as $\widetilde W(\fg)/J$, where as previously, $J$ is the maximal ideal
intersecting the local superalgebra trivially.
It is clear that the identities
(\ref{IdealJ2}) generate such an ideal. The only instance when we specialised to
$\fg=A_{n-1}$, and more specifically, the tensor structure, was in the
proof that (\ref{IdealJ2}) indeed generates the ideal $J$.
The corresponding statement for $W(D_r)$ and $W(E_r)$ remains a
conjecture.

It is straightforward to check that in $W(D_r)$ for $r\geq4$ and
$W(E_r)$ for $r\geq6$, the level $-2$ relation 
$[f_{0i},[f_1,f_{0j}]]=0$ for $i,j=3,\ldots,r$, is superfluous, although it still
generates an ideal intersecting the local part trivially. The corresponding weights are in
the same $\fg$ representation as the one of $[f_{0i},f_{0j}]$, and the
full level $-2$ part of $J$ is then generated by the relations
\begin{align}
[f_{0a},f_{0b}]&=0\,,\nn\\
[(f_{02}-f_{00}),[f_1,f_{0i}]]&=0\,,\quad i\geq3\,\label{IdealJDE}.
\end{align}

The BKM algebra with the same Dynkin diagram as $W(D_r)$ is
finite-dimensional, ${\scr B}(D_r)=\mathfrak{osp}(r,r|2)$. There are
generators at level $\ell=-2,-1,0,1,2$, with $D_r$ singlets at
$\ell=0,\pm2$. Therefore, $W(D_r)$ will not have any generators at level
$\ell\geq3$. 

We conjecture that the level decomposition of $W(D_r)$
consists of an infinitely repeating sequence of antisymmetric modules,
such that there are scalars at levels $2-2p$,
vectors at levels $1-2p$, antisymmetric 2-index tensors at levels $-2p$ etc. for $p=0,1,2,\ldots$.
The completeness of the ideal at level $-2$ should be possible to
check in this case. For $S(D_r)$, however, the scalar at level $2$
is part of the ideal $J$, and there is (by definition) no singlet at
level 0. There will be no
recurrence of the antisymmetric tensors. 
Therefore, the Lie superalgebra $S(D_r)$ is
finite-dimensional, and it is isomorphic to $H(2r)$ in the
classification of Kac \cite{Kac77B}.

In the case $\fg=E_r$,  the BKM algebra ${\scr B}(E_r)$ is
infinite-dimensional (see \eg \cite{Cederwall:2015oua}), and
$W(E_r)$ and $S(E_r)$ contain generators at all integer levels. A list of $E_r$ representations for $4 \leq r \leq 8$ up to level $12-r$
can be found for example in \cite{Howe:2015hpa}.
We conjecture that the relations (\ref{IdealJDE}) generate the maximal ideal
$J$ also in this case, but a proof is so far lacking.

\subsection{Realization of \texorpdfstring{$W(E_n)$}{WEn}}\label{WEn-section}

In Section \ref{WgSection} we assumed that the algebra $W(\fg)$ is non-trivial, \ie, that the ideal of the free Lie superalgebra $F$
generated by the relations (\ref{eigen})--(\ref{eifjf0a})
is not equal to $F$ itself, but a proper ideal. To verify the assumption it is sufficient to find a non-trivial
algebra homomorphism from $W(\fg)$ to a non-trivial algebra. When $\fg$ is finite-dimensional, it is straightforward to construct such a
homomorphism from $S(\fg)$ to the (original) tensor hierarchy algebra associated to $\fg$, given the structural details in Section \ref{WgSection}.
The homomorphism can then be extended to a homomorphism from $W(\fg)$ to the extended version of the tensor hierarchy algebra.

However, when $\fg$ is infinite-dimensional, it is not obvious that the assumption is true. The infinite-dimensional cases that we are
most interested in are $\fg=E_r$ for $r\geq 9$, with the default choice of ``node 1'', which means $\fg'=E_{r-1}$. The assumption that
both $\fg$ and $\fg'$ are simple
excludes the affine Lie algebra $E_9$ and adjusts the range of $r$ to $r \geq 11$. 
In \cite{Bossard:2017wxl} a tensor hierarchy algebra 
associated to $E_r$ was defined also for $r\geq 9$, with particular focus on the case $r=11$.
We will end this section (and the paper) by briefly giving a surjective homomorphism from 
$W(E_r)$ to an extension of the algebra defined in \cite{Bossard:2017wxl} (which is the image of the $S(E_r)$ subalgebra).

Consider the Grassmann superalgebra $\Lambda=\Lambda(n)$. Since it is an associative algebra with identity element, there is an injective
homomorphism $\Lambda(n) \to \End \Lambda(n)$ given by left multiplication. It is common to use the same notation for
the image of any element under this homomorphism as for the element itself, writing
\begin{align} \label{leftmult}
x \quad:\quad \Lambda(n) \to \Lambda(n)\,, \quad  y \mapsto xy\,.
\end{align}
We will employ this convention, but at the same time it will be important to distinguish between the two copies
of $\Lambda(n)$. Therefore, we denote the identity elements in $\Lambda(n)$ and $\End \Lambda(n)$ by $E$ and $L$, respectively,
and write out these elements explicitly in the expressions. For example, (\ref{leftmult}) then becomes 
\begin{align}
xL \quad: \quad\Lambda(n) \to \Lambda(n)\,, \quad  yE \mapsto xyE\,.
\end{align}
For any triple of indices $a,b,c=0,1,\ldots,n-1$, we define a map $F_{abc}: \Lambda \to \End \Lambda$ by
\begin{align}
F_{abc}(xE)=3(K_{[a}K_b x)K_{c]}+(-1)^{|x|}(K_{a}K_{b}K_{c}x)L,
\end{align}
where $K_b$ is the contraction
\begin{align}
K_b : \xi^{c_1} \cdots \xi^{c_q}E \mapsto q\delta_b{}^{[c_1} \xi^{c_2}\cdots \xi^{c_q]}E
\end{align}
defined in (\ref{contraction}). Set
$K^{a_1\cdots a_p}{}_b = \xi^{a_1} \cdots \xi^{a_p} K_b$ and $K=\sum_{a=0}^{n-1}K^a{}_a$.

Consider now the local Lie superalgebra
$u(\Lambda) = U_{-1} \oplus U_0 \oplus U_1$, as defined in Section \ref{universal} (but note that the $\mathbb{Z}$-grading is not consistent in this case, since the $\mathbb{Z}_2$-graded vector space $U_1=\Lambda$ is not homogeneous).
Thus
\begin{align}
U_1 &= \Lambda\,, & U_0&= \End \Lambda\,, & U_{-1} &= \mathrm{Hom}\,(\Lambda,\End \Lambda)\,.
\end{align} 
Let $w_E(n)$ be the local subalgebra of $u(\Lambda)$ generated by all elements in $U_1$, the elements
$F_{abc}$ in $U_{-1}$, and the elements $K^{a_1\cdots a_p}{}_b$ in $U_0$.
Let $W_E (n)$ be the minimal Lie superalgebra with local part $w_E(n)$. 
Consider the map $W(E_n)  \to W_E(n)$ given by
\begin{align} 
e_0 &\mapsto K_0\,,& e_i &\mapsto K^{i-1}{}_{i}\,, & f_i &\mapsto K^{i}{}_{i-1}\,,\nn
\end{align}
\begin{align}
e_n &\mapsto \xi^{(n-3)}\xi^{(n-2)}\xi^{(n-1)}E\,, & f_n &\mapsto F_{(n-3)(n-2)(n-1)}\,,\nn\\
h_0 &\mapsto K-3L-K^0{}_0\,, & h_i&\mapsto K^{i-1}{}_{i-1} -K^{i}{}_{i}\,,\nn\\
f_{00}&\mapsto K^{0} - 3L^0\,,&
f_{0i}&\mapsto K^{0(i-1)}{}_{i-1} -K^{0i}{}_{i}\,,\nn
\end{align}
\begin{align}
h_n&\mapsto K^{n-3}{}_{n-3}+K^{n-2}{}_{n-2}+K^{n-1}{}_{n-1}-L\,,\nn\\
f_{0n}&\mapsto K^{0(n-3)}{}_{n-3}+K^{0(n-2)}{}_{n-2}+K^{0(n-1)}{}_{n-1}-L^0
\end{align}
for the set of generators, where $K^0=\xi^0 K$ and $L^0=\xi^0 L$.
It is straightforward to check that this map preserves all the relations (\ref{eigen})--(\ref{eifjf0a}), and thus it is a homomorphism.

\appendix

\section{The root system of \texorpdfstring{$W(A_{n-1})=W(n)$}{W(A(n-1))}}

A weight $\lambda$ can expressed in the form
\begin{align}
\lambda=k\widetilde\Lambda_0+\sum_{i=1}^{n-1}\mu_i\widetilde\Lambda_i\,.
\end{align} 
Here $\widetilde \Lambda_0$ is a weight which is orthogonal to all simple roots $\alpha_1,\alpha_2,\ldots,\alpha_r$,
and has coefficient 1 for $\alpha_0$ when expressed in
terms of simple roots. This implies that $\widetilde \Lambda_0$ is proportional to
$\Lambda_0$. 
The inverse of the Cartan matrix $B$ given in (\ref{cartan-A-case}) is
\begin{align}
B^{-1}=\left(
\begin{matrix}
-\frac{n}{n-1}&-1&-\frac{n-2}{n-1}&\cdots&-\frac{1}{n-1}\cr
-1&0&0&\cdots&0\cr
-\frac{n-2}{n-1}&0&&&\cr
\vdots&\vdots&&{A'}^{-1}&\cr
-\frac{1}{n-1}&0&&&\cr
\end{matrix}
\right)\,,
\end{align}
where $A'$ is the Cartan matrix for $A_{n-2}$. An analogous
structure (with $0$'s in the second row and column, and the inverse
Cartan matrix for the algebra with Dynkin diagram obtained by deleting
nodes $0$ and $1$) arises also for other choices of $\fg$.
From the upper left corner we get
\begin{align}
\widetilde\Lambda_0=-\frac{n-1}{n}\lambda_0\,.
\end{align} 
This also implies that $(\widetilde\Lambda_0,\widetilde\Lambda_0)=-\frac{n-1}{n}$.
The $\widetilde\Lambda_i$'s satisfy 
$(\widetilde\Lambda_i,\alpha_j)=\delta_{ij}$, and have vanishing
$\alpha_0$ component when expressed in the basis of simple
roots. Thus, $(\widetilde\Lambda_0,\widetilde\Lambda_i)=0$.
The length of the weight $\lambda$ becomes
\begin{align}
  (\lambda,\lambda)=-\frac{n-1}{n}k^2+(\mu,\mu)\,,
  \label{lambdalength}
\end{align}
where the scalar product on the right hand side is calculated for
weights $\mu$ of $A_{n-1}$.

We can use (\ref{lambdalength}) together with the known $A_{n-1}$
  representations to give the lengths of all roots
in $W(n)$ or $S(n)$. The representation $(0\ldots010\ldots01)$, with
the first $1$ in position $k+1$, occurs at level $-k$. It contains two
Weyl orbits, represented by the dominant weights $\mu_{k+1}+\mu_{n-1}$
and $\mu_k$, where $\mu_i$ are simple $A_{n-1}$ weights. The other
  representation at level $-k$ has highest weight $\mu_k$.
  The lengths squared of these weights are
  \begin{align}
    (\mu_k,\mu_k)&=\frac{k(n-k)}{n}\,,&   
    (\mu_{k+1}+\mu_{n-1},\mu_{k+1}+\mu_{n-1})
    &=\frac{k(n-k)}{n}+2\,.
  \end{align}
  Insertion into (\ref{lambdalength}) tells us that the root
  lengths at level $-k$ are $k-k^2$ and $2+k-k^2$.
  Roots with $(\lambda,\lambda)>0$ appear at level $0$ and $-1$. 
  Null roots appear at level $1$, $-1$ and $-2$, the last case
  for $n\geq4$.

\begin{figure} 
\begin{center}
{\includegraphics[width=4in]{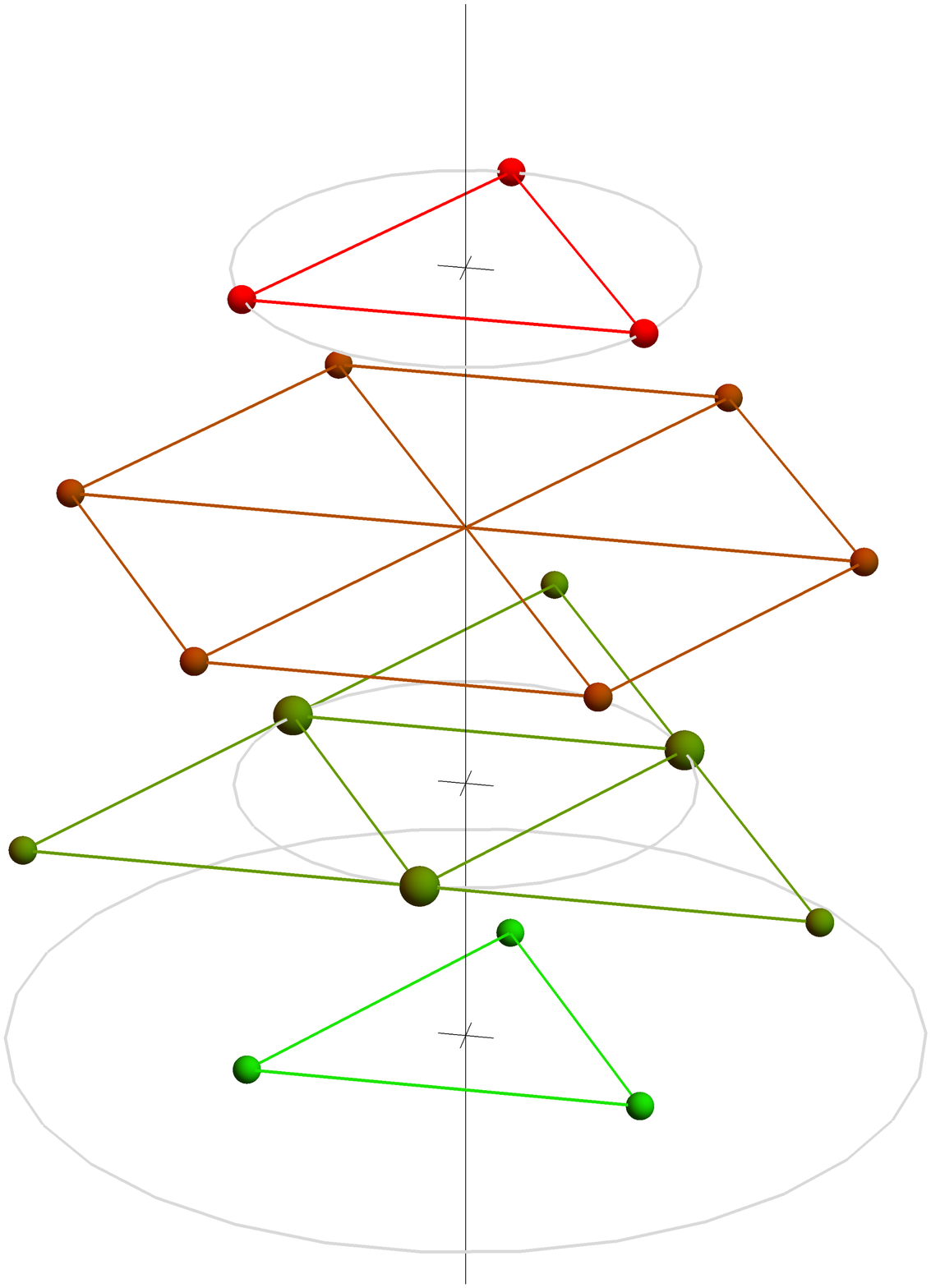}}
\caption{\label{rotated} \it The root system of $W(3)$. Each ball
  corresponds to a root, and the sizes of the balls
to the multiplicities of the
roots. The roots at each 
plane are the weights of the $A_2$-representation that occurs at the
corresponding level in $W(3)$.
The circles are the intersections of the planes at the different
levels with the ``light cone'' consisting of null weights.} 
\end{center}
\end{figure}

\begin{figure}
\begin{center}
{\includegraphics[width=1.48in]{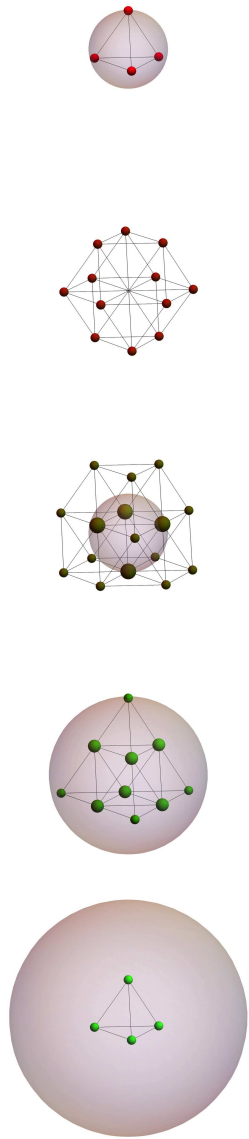}}
\caption{\label{w4all}\it The root system of $W(4)$, divided into
  levels, from 1 at the top to $-3$ at the bottom. The spheres
  indicate the intersections of the level planes with the
  cone of null weights. Note the presence of null roots
  at levels $1$, $-1$ and $-2$.
}
\end{center}
\end{figure}

The root system for $W(3)$ is depicted in Figure \ref{rotated},
and the one for $W(4)$ in Figure \ref{w4all}.
The $W(3)$ roots are listed in Table \ref{w3table}.

\bibliographystyle{utphysmod2}

\providecommand{\href}[2]{#2}\begingroup\raggedright\endgroup

\newpage

\begin{center}
  \begin{table}[h]   
\begin{tabular}{| c  | c |  c | c  | c |c|}	

\hline
level &   basis  & $\mathfrak{sl}(3)$ representation & roots $\alpha$ & mult $\alpha$ & $(\alpha,\alpha)$\\ 
\hline
$1$ &  $\begin{matrix} K_0 \\ K_1\\K_2 \end{matrix}$  & $\overline{\mathbf{3}}$ &
$\begin{matrix}
\alpha_0\\
\alpha_0+\alpha_1\\
\alpha_0+\alpha_1+\alpha_2
\end{matrix}$
&$\begin{matrix}1\\1\\1\end{matrix}$&$\begin{matrix}0\\0\\0\end{matrix}$\\ 
\hline

$0$ &  $\begin{matrix} K^0{}_1 \\ K^1{}_2\\K^0{}_2\\K^1{}_0\\K^2{}_1\\K^2{}_0 \end{matrix}$  & ${\mathbf{8}} +{\mathbf{1}}$ &$\begin{matrix}
\alpha_1\\
\alpha_2\\
\alpha_1+\alpha_2\\
-\alpha_1\\
-\alpha_2\\
-\alpha_1-\alpha_2
\end{matrix}$&$\begin{matrix}1\\1\\1\\1\\1\\1\end{matrix}$&$\begin{matrix}2\\2\\2\\2\\2\\2\end{matrix}$\\ 
\hline
$-1$ &  $\begin{matrix} K^{01}{}_2 \\ K^{20}{}_1\\K^{12}{}_0\\K^{01}{}_1 \ \ \, K^{02}{}_2 \\K^{12}{}_2 \ \ \, K^{10}{}_0\\K^{20}{}_0 \ \ \, K^{21}{}_1 \end{matrix}$  &
$\overline{\mathbf{6}}+{\mathbf{3}}$ &$\begin{matrix}
-\alpha_0+\alpha_2\\
-\alpha_0-\alpha_2\\
-\alpha_0-2\alpha_1-\alpha_2\\
-\alpha_0\\
-\alpha_0-\alpha_1\\
-\alpha_0-\alpha_1-\alpha_2
\end{matrix}$&$\begin{matrix}1\\1\\1\\2\\2\\2\end{matrix}$&$\begin{matrix}2\\2\\2\\0\\0\\0\end{matrix}$\\  
\hline
$-2$ &    $\begin{matrix} K^{012}{}_2 \\ K^{012}{}_1\\K^{012}{}_0\end{matrix}$ & $\overline{\mathbf{3}}$ &$\begin{matrix}
-2\alpha_0-\alpha_1\\
-2\alpha_0-\alpha_1-\alpha_2\\
-2\alpha_0-2\alpha_1-\alpha_2\\
\end{matrix}$&$\begin{matrix}1\\1\\1\end{matrix}$&$\begin{matrix}-2\\-2\\-2\end{matrix}$\\ 
    \hline
\end{tabular}
\vspace{12pt}
\caption{\it The $W(3)$ root system.\label{w3table}}
\end{table}
\end{center}

\end{document}